\documentclass{amsart}
\usepackage{amssymb} 
\usepackage[all]{xy} 
\newcommand\Hd{\operatorname{HH}} 
\newcommand\GGg[1]{{\mathbb G}(#1)} 
\newcommand\GG{{\mathbb G}}  
\newcommand\cohom{\operatorname{H}}  
\newcommand\mfk{\mathfrak}  
\newcommand\rg{\operatorname{reg}}  
\newcommand\nil[1]{{\mathcal U}_{#1}}  
\newcommand\class{\langle S \rangle} 
\newcommand\ind{\operatorname{ind}} 
\newcommand\rk{{\,\operatorname{rk}}} 
\newcommand\cunip[2]{\cohom_{#1} (#2,\CI_c({#2}_u)_{\delta})} 
\newcommand\orb{\mathcal O}

\newcommand\datver[1]{\def\datverp% 
 {\par\boxed{\boxed{\text{Version: #1; Run: \today}}}}} 
\datver{3.0A; Revised: 8/16/2000}

\newcommand\CC{\mathbb C} 
\newcommand\FF{\mathbb F} 
\newcommand\KK{\mathbb K} 
\newcommand\NN{\mathbb N} 
\newcommand\QQ{\mathbb Q} 
\newcommand\RR{\mathbb R} 
\newcommand\ZZ{\mathbb Z}

\newcommand\CI{{\mathcal C}^{\infty}} 
\newcommand\CIc{{\mathcal C}^{\infty}_{\text{c}}} 
\newcommand\Mand{\text{ and }} 
\newcommand\ie{i.e.,\ } 
\newcommand\inff{\operatorname{inf}} 

\newtheorem{theorem}{Theorem}

\newtheorem{proposition}{Proposition} 
\newtheorem{corollary}{Corollary} 
\newtheorem{lemma}{Lemma} 
 
\theoremstyle{definition} 
\newtheorem{definition}{Definition} 
 
\theoremstyle{remark} 
\newtheorem{remark}{Remark} 
\begin{document} 
 
\title[Homology of $p$-adic groups] 
{Higher orbital integrals, Shalika germs, and the Hochschild
homology of Hecke algebras} 

\author[V. Nistor]{Victor Nistor$^{1}$} 
	\address{Pennsylvania
        State University, Math. Dept., University Park, PA 16802, USA}
	\email{nistor@math.psu.edu}
        \thanks{$^1$partially supported by a NSF Young Investigator
        Award DMS-9457859 and a Sloan Research Fellowship. Manuscripts
        are available from {\bf
    http:{\scriptsize//}www.math.psu.edu{\scriptsize/}nistor{\scriptsize/}}}

\dedicatory\datverp 
%\date\datverp 
 
\begin{abstract}  
We give a detailed calculation of the Hochschild and cyclic homology
of the algebra $\CIc(G)$ of locally constant, compactly supported
functions on a reductive $p$--adic group $G$.  We use these calculations to
extend to arbitrary elements the definition the higher orbital
integrals introduced by Blanc and Brylinski for regular semisimple
elements. Then we extend to higher orbital integrals some results of
Shalika. We also investigate the effect of the ``induction morphism''
on Hochschild homology.
\end{abstract} 

\maketitle \tableofcontents 
 
\section*{Introduction} 
 
Orbital integrals play an important role 
in the harmonic analysis of a reductive $p$--adic group $G$;
they are, for instance, one of the main ingredients in the Arthur--Selberg
trace formula.  Orbital integrals on unimodular groups are
a particular case of invariant distribution, which have been
used in \cite{Bernstein} to prove the irreducibility of 
certain induced representations of $GL_n$ over a $p$-adic
field. 
 
By definition, an invariant distribution on a unimodular group $G$
gives rise to a trace (\ie a Hochschild cocycle of degree zero) on
$\CIc(G)$, the Hecke algebra of compactly supported, locally constant,
complex valued functions on $G$.  It is interesting then to try to
analyse as completely as possible the Hochschild homology and
cohomology groups of the algebra $\CIc(G)$ (denoted $\Hd_*(\CIc(G))$
and $\Hd^*(\CIc(G))$, respectively). 

In this paper, $G$ will be the set of $\FF$-rational points of a
linear algebraic group $\GG$ defined over a finite extension 
$\FF$ of the field $\QQ_p$ of $p$-adic numbers, $p$ being a 
fixed prime number. The group $\GG$ does not have to be reductive,
although this is certainly the most interesting case. When we shall
assume $\GG$ (or $G$, by abuse of language) to be reductive, we shall
state this explicitely.  
For us, the most important topology to consider on 
$G$ will be the locally compact topology induced from
an embedding of $G \subset GL_n(\FF)$. Nevertheless, the Zariski 
topology on $G$ will also play a role in our study.

To state the main result of this paper, we need to introduce first the
concept of standard subgroup. For any set $A \subset G$, we shall denote 
\begin{equation*}
	C(A) := \{ g \in G, ga=ag, \ \forall a \in A \}
\end{equation*}
and $Z(A) := A \cap C(A)$. This latter notation will be used only when
$A$ is a subgroup of $G$. A commutative subgroup $S$ of $G$ is
called {\em standard} if $S = Z(C(s))$ for some semi-simple element 
$s \in G$. Our results will be stated in terms of standard subgroups of 
$G$. We shall denote by $H_u$ the set of unipotent elements of a subgroup $H$.
Sometimes, the set $C(S)_u$ is also denoted by $\nil S$, in order to avoid
having to many paranthesis in our formulae.  One of the
main results of this paper (Theorem~\ref{Theorem.Gen}) identifies the
groups $\Hd_*(\CIc(G))$ in terms of the following
data: the set $\Sigma$ of (conjugacy classes of) standard subgroups
$S$ of $G$, the subset $S^{\rg} \subset S$ of $S$-regular elements, the
action of the Weyl group $W(S)$ of $S$ on $\CIc(S)$, and the
continuous cohomology of the $C(S)$--module
\begin{equation*}
	\CIc(\nil S)_{\delta}:=\CIc(C(S)_u) \otimes \Delta_{C(S)},
\end{equation*} 
where $\Delta_{C(S)}$ denotes the modular function of the group 
$C(S)$. More precisely, if $G$ is a $p$-adic group defined over a field of 
characteristic zero, as before, then Theorem \ref{Theorem.Gen} states the existence
of an isomorphism
\begin{equation} \label{eq.result}
	\Hd_q(\CIc(G)) \simeq \bigoplus\limits_{S \in \Sigma}
	\CIc(S^{\rg})^{W(S)} \otimes \cohom_q(C(S), \CIc(\nil S)_\delta),
\end{equation}
which can be made natural by using a generalization of the Shalika germs.

It is important to relate this result with the periodic cyclic
homology groups of $\CIc(G)$. For the Hecke algebra $\CIc(G)$, the
periodic cyclic homology is related to Hochschild homology by
\begin{equation}
	{\rm HP}_q(\CIc(G)) \simeq \oplus_{k \in \ZZ} \Hd_{q +
	2k}(\CIc(G))_{comp},
\end{equation}
that is, ${\rm HP}_*(\CIc(G))$ is the localization of Hochschild
homology to the $G$-invariant subset of compact elements of
$G$. This relation is implicit in \cite{HN}.
Consequently, the results of this paper complement the results on
the cyclic homology of $p$--adic groups in \cite{HN,Schneider}. It is
interesting to remark that ${\rm HP}_*(\CIc(G))$ can also be described
in terms of the admissible spectrum of $G$, see \cite{KNS}, and hence
our results have significance for the representation theory of $p$-adic 
groups. See also \cite{N-algg} for similar results on (the groups of real 
points of) algebraic groups defined over $\RR$. These periodic cyclic
cohomology groups are isomorphic to $K_*(C_r^*(G))$, by combining results 
from \cite{BHP}, \cite{Lafforgue}, and \cite{HN}.
   
Assume for the moment that $G$ is reductive. Then, in order to better understand 
the role played by the groups $\Hd_*(\CIc(G))$ and $\cohom_*(G,\CIc(G_u))$ in the 
representation theory of $G$, we relate $\cohom_*(G,\CIc(G_u))$ to the analogous 
cohomology groups, $\cohom_*(P,\CIc(P_u)_{\delta})$ and $\cohom_*(M,\CIc(M_u)_{\delta})$,  
associated to parabolic subgroups $P$ of $G$ and to their Levi 
components $M$. In particular, we define morphisms between these Hochschild
homology groups that are analogous to the 
induction and inflation morphisms that play such a prominent role 
in the representation theory of $p$--adic groups. These morphisms are induced 
by morphisms of algebras.  
 
In \cite{Blanc-Brylinski}, Blanc and Brylinski have introduced higher
orbital integrals associated to regular semisimple elements by proving
first that
\begin{equation}\label{eq.MacLane}
	\Hd_q(\CIc(G)) \simeq \cohom_q(G, \CIc(G)_{\delta}), 
\end{equation}
a result which they dubed ``the MacLane isomorphism.''
(Actually, they did not have to twist with the modular function, because
they worked only with unimodular groups $G$, see Lemma \ref{Lemma.BlBr}
for the slightly more general version needed in this paper), after that 
we rely more on filtrations of the $G$-module
$\CIc(G)$, rather than on localization. This allows us to define
higher orbital integrals at arbitrary elements. Then, we study the
properties of these orbital integrals and we obtain in particular a
proof of the existence of abstract Shalika germs for the higher
orbital integrals. Actually, the existence of Shalika germs turns out
to be a consequence of some general homological properties of the ring
of invariant, locally constant functions on the group $G$. We also use 
the techniques developed in \cite{N-algg} in the framework of real algebraic 
groups. It would be interesting to relate the results of this paper to those 
of \cite{Kazhdan-cusp}.

This is the revised version of an preprint that was first circulated in February 
1999. I would like to thank Paul Baum, David Kazhdan, Robert Langlands, George Lusztig, 
Roger Plymen, and Peter Schneider for useful comments and discussions.

\section{Homology of Hecke algebras\label{Sec.Homology}}

In this section we shall use several general results on Hochschild homology 
of algebras, on algebraic groups, and on the continuous cohomology  
of totally disconnected groups. Good references are  
\cite{Borel1,Borel-Wallach,MacLane1}, for the general theory, and \cite{KNS}
for questions related to Hochschild homology. 
 
If $G$ is a group and $A\subset G$ is a subset, we denote by  
$C(A)$ the {\em centralizer} of $A$, that is,   
the set of elements of $G$ that commute with every element of $A$, and by   
$N(A)$ the {\em normalizer} of $A$, that is, the set of elements $g\in G$   
such that $gAg^{-1} =A$. By $Z=Z(G)=C(G)$ we denote the center of $G$. 
 
If $X$ is a totally disconnected, locally compact, non-discrete space $X$, we   
denote by $\CIc(X)$ the space of compactly supported, locally constant,
complex valued  functions on $X$. Recall that, if $U \subset X$ is an 
open subset of $X$ as above, then restriction defines an isomorphism  
\begin{equation}\label{eq.disc.exact} 
\CIc(X) / \CIc(U) \simeq \CIc(X \setminus U). 
\end{equation}    
  
Let $\GG$ be a linear algebraic group defined over a 
totally disconnected locally compact field $\FF$. Thus $\FF$ is 
a finite algebraic extension of  $\QQ_p$, the field 
of $p$-adic numbers. The set $\GGg \FF$ of  
$\FF$-rational points of $\GG$ is called a $p$--adic group 
and will be denoted simply by $G$.   
It is known \cite{Borel1} that $G=\GGg \FF$ identifies with a closed  
subgroup of $GL_n(\FF)$, and hence it has a natural locally compact  
topology that makes it a totally disconnected space. We fix a Haar  
measure $dg$ on $G$.      
 
Consider now the space $\CIc(G)$ of compactly supported, locally constant 
functions on $G$. Fix a Haar measure $dh$ on $G$. Then the 
convolution product, denoted $*$, is defined by    
\begin{equation*}  
	f_1 * f_2 (g) =\int_G f_1(h) f_2(h^{-1}g)dh,  
\end{equation*}  
makes $\CIc(G)$ an algebra, the {\em Hecke algebra} of $G$. It is  
important in representation theory to determine the ($Ad_G$--)invariant 
linear functionals on $\CIc(G)$. If $G$ is unimodular, the space of invariant 
linear functionals on $\CIc(G)$
coincides with the space of traces on $\CIc(G)$; moreover, since the space
of traces of $\CIc(G)$ identifies with $\Hd^0(\CIc(G))$, the first Hochschild 
cohomology group of $\CIc(G)$, it is reasonable to ask what are
the groups $\Hd^q(\CIc(G))$, the Hochschild homology groups of $\CIc(G)$, 
in general. Since $\Hd^q(\CIc(G))$ is the algebraic
dual of $\Hd_q(\CIc(G))$, it is enough to concentrate on  Hochschild
homology. The computation of the groups $\Hd_q(\CIc(G))$ is the main
purpose of this paper. 

We first recall the definition of $\Hd_*(\CIc(G))$.
Let 
$$
	\CIc(G^{q+1}) =\CIc(G) \otimes \CIc(G) \otimes \ldots \otimes 
	\CIc(G),
$$ 
$(q+1)$--times, be the usual (algebraic) tensor product of 
vector spaces. The Hochschild differential $b:\CIc(G^{q+2}) \to 
\CIc(G^{q+1})$ is  given by  
\begin{multline} 
\label{eq.Hochschild} 
	(bf)(g_0,g_1,\ldots,g_q)= \sum_{j=0}^q (-1)^j \int_G 
	f(g_0,\ldots, g_{j-1}, \gamma, \gamma^{-1}g_j, 
	g_{j+1}, \ldots, g_q) d\gamma \\ 
	+(-1)^{q+1}\int_G f(\gamma^{-1}g_0, g_1, \ldots, g_q, \gamma) d\gamma. 
\end{multline} 
By definition, {\em the $q$-th Hochschild homology group} of $\CIc(G)$,
denoted $\Hd_q(\CIc(G))$, is the $q$-th homology group of the complex
$(\CIc(G^{q+1}),b)$.  Since $\CIc(G)$ is an inductive limit of unital
algebras, this definition coincides with the usual definition of
Hochschild homology for non-unital algebras (using algebras with adjoint
unit).
 
The group $G$ acts by conjugation on $\CIc(G)$, and we denote by   
$\CIc(G)_{\operatorname{ad}}$ the $G$-module defined by this action. Also,
let $\Delta_G$ denote the modular function of $G$, which we recall, is
defined by the relation
$$
	\Delta_G(h)\int_G f(gh)dg = \int_G f(g)dg.
$$
We shall be especially interested in the $G$--module 
$\CIc(G)_{\delta}$ obtained from $\CIc(G)_{\operatorname{ad}} $ by 
twisting it with the modular function.
More precisely, let $\CIc(G)_{\delta} = \CIc(G)$
as vector spaces, and let the action of $G$ on functions be given
by the formula
\begin{equation}
	(\gamma \cdot f)(g)= \Delta_G(\gamma) f(\gamma^{-1} g \gamma), 
	\qquad f \in \CIc(G)_{\delta}.
\end{equation}
The reason for this twisting is that, for $G$ non-unimodular, the
traces of $\CIc(G)$ are the $G$--invariant functionals on $\CIc(G)_{\delta}$,
{\em not} on $\CIc(G)$ (this is an immediate consequence of Lemma 
\ref{Lemma.BlBr}). More generally, our approach to the Hochschild homology of 
$\CIc(G)$ is based on Lemma \ref{Lemma.BlBr}. 

Before stating and proving Lemma \ref{Lemma.BlBr}, we need to
introduce some notation. First, if $M$ is an arbitrary $G$--module, we
denote by $M \otimes \Delta_G$ the tensor product of the $G$--modules
$M$ and $\CC$, where the action on $\CC$ is given by the multiplication
with the modular function of $G$. (We sometimes write this as
$\CIc(G)_{\delta} = \CIc(G) \otimes \Delta_G$.) 

If $M$ is a right $G$--module and $M'$ is a left
$G$--module, then $M \otimes_G M'$ is defined as the quotient of
$M \otimes M'$ by the submodule generated by $mg \otimes m'-
m \otimes gm'$. For example, if $H \subset G$ is a closed subgroup
and  if $X$ is a left $H$--space, then we have an isomorphism of
$G$--spaces
\begin{equation} \label{eq.ind.quot}
\CIc(G) \otimes_H (\CIc(X) \otimes \Delta_H) \simeq \CIc(G \times_H X),
\end{equation}
where $G \times X$ is the quotient $(G \times X)/H$ for the action 
$h(g,x)=(gh^{-1},hx)$. This isomorphism is obtained by observing that 
the natural map
\begin{multline*}
t_X : \CIc(G) \otimes \CIc(X) = \CIc(G \times X) \longrightarrow
\CIc(G \times_H X),\\
t_X(f)(\overline{(g,x)}) = \int_H f(gh,h^{-1}x) dh,
\end{multline*} 
passes to the quotient to give the desired isomorphism. Sometimes
it will be convenient to regard a left $G$--module as a 
right $G$--module by replacing $g$ with $g^{-1}$.

Also, recall that a $G$--module $M$ is {\em smooth} if, and only if, 
the stabilizer of each element of $M$ is open in $G$.
Then one can define the {\em continuous homology} groups of $G$
with coefficients in the smooth module $M$, denoted 
$\cohom_k(G,M)$, using tensor products as follows. Let
${\mathcal B}_q(G)=\CIc(G^{q+1})$, $q=0,1,\ldots$,
be the Bar complex of the group $G$, 
with differential 
$$ 
	(df)(g_0,g_1,\ldots,g_q) = \sum_{j =0}^{q+1}(-1)^j\int_G 
	f(g_0,\ldots,g_{j-1},\gamma,g_{j},\ldots,g_q)d\gamma. 
$$ 
Then the complex $({\mathcal B}_q,d)$ gives a resolution of $\CC$ with 
projective $\CIc(G)$--modules, and the complex  
\begin{equation} \label{eq.cont.cohom}
{\mathcal B}_q(G)\otimes_G  M
\end{equation}  
computes $\cohom_q(G,M)$. See \cite{Blanc,Borel-Wallach}. 

We shall need the following extension of a result from \cite{Blanc-Brylinski}:

\begin{lemma}\label{Lemma.BlBr}
Let $\CIc(G)_{\delta}=\CIc(G) \otimes 
\Delta_G$ be the $G$--module obtained by twisting the adjoint action
of $G$ on $\CIc(G)$ by the modular function. Then we have a natural isomorphism
\begin{equation}\label{eq.BlBr} 
\Hd_q(\CIc(G)) \simeq \cohom_q(G,\CIc(G)_{\delta}). 
\end{equation} 
\end{lemma}

\begin{proof} 
Consider the complex \eqref{eq.cont.cohom}, which computes the
continuous cohomology of $M =\CIc(G)_{\delta}$, and
let $h_G:{\mathcal B}_q(G) \otimes\CIc(G)_{\delta}\simeq 
\CIc(G) \otimes \CIc(G^{q+1}) =\CIc(G^{q+2}) \to \CIc(G)^{\otimes q+1}$ 
be the map 
\begin{multline*}
	h_G(f)(g_0,g_1,\ldots,g_q) 
	=\int_G f(g^{-1}h g ,\,g^{-1}g_0,\,g^{-1}g_0g_1,\,\ldots,\, 
	g^{-1}g_0g_1\ldots g_q) dg,\\ h = g_0g_1 \ldots g_q .
\end{multline*}

As in equation \eqref{eq.ind.quot}, the map $h_G$ descends to the quotient to
induce an isomorphism
\begin{equation} \label{eq.tensor} 
	\tilde{h}_G:{\mathcal B}_q(G) \otimes_G \CIc(G)_{\delta}\simeq
	\CIc(G^{q+2}) \otimes_G \CC \simeq \CIc(G)^{\otimes q+1}
\end{equation} 
of complexes, that is $\tilde{h}_G \circ (d \otimes_G 1) =b \circ
\tilde{h}_G$, which establishes the isomorphism
$\cohom_q(G,\CIc(G)_{\delta}) \simeq \Hd_q(\CIc(G))$, as 
desired.
\end{proof}

To better justify the twisting of the module $\CIc(G)$ by the modular
function in the above lemma, note that the trivial representation of
$G$ gives rise to an obvious morphism $\pi_0:\CIc(G) \to \CC$, by
$\pi_0(f)=\int_G f(g)dg$, which hence defines a trace on $\CIc(G)$.
However $\pi_0$ is not $G$--invariant for the usual action of $G$, but
is invariant if we twist the adjoint action of $G$ by the modular
function, as indicated.
  
We proceed now to a detailed study of the $G$-module
$\CIc(G)_{\delta}$.  
First we define a natural $Ad_G$-invariant stratification of
$G$, called the {\em standard stratification} of $G$.
  
Let $\mfk g$ be the Lie algebra of $G$ in the sense of linear
algebraic groups.  Denote by $a_i(g)$ the coefficients of the
polynomial $\operatorname{det}(t +1 - Ad_g)$,
$$  
	\operatorname{det}(t +1 - Ad_g)= \sum_{i=0}^{m} a_i(g) t^i \in
	\FF[t].
$$  
Let $a_r$ be the first non-zero coefficient $a_i$, and define   
$$ 
	V_k=\{a_r=a_{r+1}=\ldots = a_{r+k-1}=0\}.
$$  
Thus $V_0=G$, by convention, and  
$G \setminus V_1=G'$, the set of regular elements of $G$ 
if $G$ is reductive. Also, $V_{m+1} = \emptyset$ because $a_{m}=1$.
We observe that the functions $a_i(g)$  
are $G$-invariant polynomial functions on $G$, and that they depend only  
on the semisimple part of $g$.  
 
The description of the Hochschild homology of $\CIc(G)$ that we shall
obtain is formulated in terms of certain commutative subgroups 
of $G$, called {\em standard}, that we now define.

\begin{definition} \label{Def.standard}
A commutative subgroup $S \subset G$ is called  
{\em standard} if, and only if, there exists a semisimple element $s_0 \in G$ 
such that $S$ is the center of $C(s_0)$, the centralizer of $s_0$ in $G$.  
\end{definition}

If $s_0$ and $S$ are as in the definition above, then $S$ is the set 
of $\FF$--rational points of a subgroup ${\mathbb S}$ of $\GG$ and 
$C(S)=C(s_0)$. The element $s_0$ in the above definition is not 
unique in general, and, for a standard subgroup $S \subset G$,   
we denote by $S^{\rg} \subset S$ the set of semisimple elements   
$s \in S$ such that $C(S)=C(s)$. An element $s \in S^{\rg}$ will
be called {\em $S$-regular}. This set is not  
empty by the definition of a standard subgroup.  Note that $S^{\rg}$ 
depends on $G$ also, and not only on $S$. 
 
Standard subgroups exist. Indeed, if $\gamma \in G$ is semisimple, then  
$S(\gamma) := Z(C(\gamma))$, the center of the centralizer of $\gamma$,  
is a standard   
subgroup of $G$. In particular, every semisimple element of $G$ belongs  
to $S^{\rg}$, for some standard subgroup $S$ of $G$.  
For all standard subgroups $S$, the set $S^{\rg}$ is an open   
subset of $S$  in the Zariski topology.  
 
For any $p$--adic group $H$, we denote by $H_u$ the set of unipotent 
elements of $H$, and call it the {\em unipotent variety} of $H$. 
In the particular case of $H=C(S)$, where $S \subset G$ is a standard  
subgroup, we also denote  
$C(S)_u =\nil S$.   
  
In order to proceed further, recall that the Jordan decomposition of
an element $g \in G$ is $g = g_sg_u$, where $g_s$ is semisimple, $g_u$
is unipotent, and $g_sg_u=g_ug_s$. This decomposition is unique
\cite{Borel1}.  If $g=g_sg_u$ is the Jordan decomposition of $g \in G$
and if $g_s \in S^{\rg}$, then $g_u \in \nil S$, by definition, and hence
$g \in S^{\rg} \nil S$.

Consider now a standard subgroup $S \subset G$ and let  
\begin{gather*}  
F_S = Ad_G( S^{reg}), \Mand \\  
F_S^u = Ad_G( S^{reg} \nil S)  
\end{gather*}  
be the set of semisimple elements of $G$ conjugated to an element of  
$S^{\rg}$ and, respectively, the set of elements $g \in G$ 
conjugated to an element of $S^{\rg}\nil S$ (\ie such that the
semisimple part of $g$ is in $F_S$).  Also, let $N(S)$ be the normalizer 
of $S$ and $W(S)=N(S)/C(S)$. Since $N(S)$ leaves $S^{\rg}$ invariant 
and is actually the normalizer of this set, it follows that the 
quotient $W(S)$ can be identified with a set of automorphisms 
of $S_{ss}$, the subgroup of semisimple elements of $S$. Since $N(S)$ 
is the set of $\FF$--rational points of an algebraic group, the 
rigidity of tori (\cite{Borel1}, page 117) shows that $W(S)$ is finite.  
 
The natural map $(g,s) \to gsg^{-1}$ descends to a map  
$$\phi_S:(G \times S^{reg})/N(S) \ni (g,s) \to gsg^{-1} \in F_S.$$ 
Similarly, we obtain a map  
$$\phi_S^u: (G \times S^{reg}  
\nil S)/N(S) \ni (g,su) \to gsug^{-1} \in F_S^u.$$  
 
In the following proposition we consider the locally compact (and  
Hausdorff) topology of $G$, and {\em not} the Zariski topology. 
Denote by $G_{ss}$ the set of semisimple elements of $G$. 
%For $G$ reductive
%and connected (in the sense of algebraic groups), most
%of the following proposition follows from results in \cite{Vigneras} (note
%however that there exists a gap in the proof of Lemma 2.7, page 956, which
%claims that the Weyl group associated to a non-semisimple element is
%finite).

\begin{proposition} \label{Prop.stratification} 
Let $S$ be a standard subgroup of $G$. Using the  
above notation, we have:   
  
(i) The set $F_S$ is an analytic submanifold of $G$, 
and the maps $\phi_S$ and $\phi_S^u$ are homeomorphisms.    
   
(ii) For each $k$, the set $V_k \setminus V_{k+1}$ is a disjoint union
of sets of the form $F_S^u$, and each $F_S^u \subset V_k \setminus
V_{k+1}$ is an open subset of $V_k$. Moreover, the set $G_{ss} \cap
(V_k \setminus V_{k+1})$ is a disjoint union of sets of the form
$F_S$, and each $F_S \subset V_k \setminus V_{k+1}$ is an open subset
of $G_{ss} \cap V_k$.
\end{proposition}

\begin{proof} (i) First we check that  $\phi_S$ and $\phi_S^u$ are  
injective. Indeed, assume that $g_1s_1 g_1^{-1}= g_2s_2 g_2^{-1}$, for  
some $s_1, s_2 \in S^{\rg}$. Then, if $g=g_2^{-1}g_1$,  
we have  
$$ 
gC(s_1)g^{-1}=C(s_2) \Rightarrow gC(S)g^{-1}=C(S)  
\Rightarrow gSg^{-1}=S, 
$$ 
and hence $g \in N(S)$. Consequently, we have
$(g_1,s_1)=(g_2g,g^{-1}s_2g)=g^{-1}(g_2,s_2)$, with $g \in N(S)$, as
desired. The same argument shows that if $F_S$ and $F_{S'}$ have a
point in common, then the standard subgroups $S$ and $S'$ are
conjugated in $G$.
 
The injectivity of $\phi_S^u$ follows from the injectivity of
$\phi_S$, indeed, if $g_1(s_1u_1) g_1^{-1}= g_2(s_2u_2) g_2^{-1}$, let
$g=g_2^{-1}g_1$ as above, and conclude that $gs_1g^{-1}=s_2$, by the
uniqueness of the Jordan decomposition. As above, this implies that
$g\in N(S)$.
 
Since the differential $d\phi_S$ is a linear isomorphism onto its
image (i.e. it is injective)
and $\phi_S$ is injective, it follows that $\phi_S$ is a local
homeomorphism onto its image (for the locally compact topologies), and
that its image is an analytic submanifold (\cite{Schikhof}, p. 38, Theorem (2.3)).  
The set $G_{ss} \cap (V_k \setminus V_{k+1})$ is an algebraic variety on which
$G$ acts with orbits of the same dimension, and hence $\phi_S$ is
proper. This proves that $\phi_S$ is a homeomorphism. Using an inverse for 
$\phi_S$, we obtain that $\phi_S^u$ is also a homeomorphism.
 
To prove now (ii), consider a standard subgroup $S \subset G$, and let
$d$ be the dimension of $C(S)$. Then $a_0=a_1=\ldots=a_{d-1}=0$ on
$S$, and $S\cap \{a_d\not = 0\}$ contains $S^{\rg}$ as an open
component. It follows that, if $s \in G_{ss} \cap (V_k \setminus
V_{k+1}) \cap S^{\rg}$, then $F_S \subset V_k \setminus V_{k+1}$. This
shows that $G_{ss} \cap (V_k \setminus V_{k+1})$ is a union of sets of
the form $F_S$. This must then be a disjoint union because the sets
$F_S$ are either equal or disjoint, as proved above.
 
Now, if $g \in V_k \setminus V_{k+1}$ has semisimple part $s$, then $s
\in F_S \subset G_s \cap (V_k \setminus V_{k+1})$, for some standard
subgroup $S$, and hence $g \in F_S^u \subset V_k \setminus
V_{k+1}$. The sets $F_S^u$ are open in the induced topology because
the map $V_k \setminus V_{k+1} \to G_s \cap (V_k \setminus V_{k+1})$
is continuous.
\end{proof}

See also \cite{Vigneras}.
  
Let $R^{\infty}:=R^{\infty}(G)$ be the ring of locally constant  
$Ad_G$-invariant functions on $G$ with the pointwise product, which we regard  
as a subset of the set of endomorphisms of the $G$--module   
$\CIc(G)_{\delta}=\CIc(G) \otimes \Delta_G$. For each $k \geq 1$, denote by   
$I_k \subset R^{\infty}$ the ideal generated by functions  
$f: G \to \CC$ of the form  
$$ 
	f=\phi(a_r,a_{r+1},\ldots,a_{r+k-1}),
$$ 
where $\phi$  
is a locally constant function $\phi: \FF^k \to \CC$ such that   
$\phi(0,0,\ldots,0)=0$. (Recall that each of the polynomials  
$a_0,\ldots,a_{r-1}$ is the $0$ polynomial.)   By convention, we set 
$I_0 = (0)$; also, it follows that $I_{m+1}=R^{\infty}$. 

Fix now $k$, and let $p_n \in I_k$ be the function   
$\phi_n(a_r,a_{r+1},\ldots,a_{r+k-1})$, where    
$\phi_n: \FF^k \to \CC$ is $1$ on the set   
$$ 
\{\xi=(\xi_0,\ldots,\xi_{k-1})\in \FF^k,\, \max |\xi_i| \geq q^{-n}\} 
$$  
and vanishes outside this set. (Here $q$ is the number of elements of  
the residual field of $\FF$, and the non-archimedean norm ``$| \:\, |$'' is   
normalized such that its range is $\{0\} \cup \{ q^n, n \in \ZZ\}$.)  
Then $p_n =p_n^2 = p_n p_{n+1}$ and $I_k = \cup p_nR^{\infty}$.  
For further reference, we state as a lemma a basic property 
of the constructions we have introduced. 
  
\begin{lemma} \label{Lemma.exact} \ 
If $M$ is a $R^\infty$--module, then $I_kM=\cup p_n M$. 
\end{lemma} 
\vspace*{2mm}  
  
As a consequence of this above lemma, we obtain the following  
result.

\begin{corollary} \label{Cor.sum} 
Consider the $G$--module 
$$M_k= (I_{k+1}/I_k) \otimes_{R^{\infty}}\CIc(G)_{\delta} \simeq 
I_{k+1}\CIc(G)_{\delta}/I_{k}\CIc(G)_{\delta}.
$$ 
Then, for each $q \ge 0$, we have an isomorphism  
$$  
\cohom_q(G,\CIc(G)_{\delta}) \simeq \bigoplus\limits_{k=0}^{m}   
\cohom_q(G, M_k)  
$$  
of vector spaces.  
\end{corollary}

\begin{proof}  
There exists a (not natural) isomorphism 
$$  
	\cohom_q(G,\CIc(G)_{\delta}) \simeq \bigoplus\limits_{k=0}^{m}
	I_{k+1} \cohom_q(G,\CIc(G)_{\delta})/I_{k}
	\cohom_q(G,\CIc(G)_{\delta}),
$$  
of vector spaces.

By the above lemma, the inclusion of $I_k \CIc(G)_{\delta}
\to\CIc(G)_{\delta}$ of $G$-modules induces natural isomorphisms
\begin{multline*} 
\cohom_q(G, I_k\CIc(G)_{\delta}) \simeq  
\cohom_q(G,\displaystyle{\lim_{\to}}\;p_n\CIc(G)_{\delta}) \\ 
\simeq \displaystyle{\lim_{\to}}\; p_n\cohom_q(G,\CIc(G)_{\delta}) 
\simeq I_k \cohom_q(G,\CIc(G)_{\delta}), 
\end{multline*} 
because the functor $\cohom_q$ is compatible with inductive limits and
with direct sums.
 
The naturality of these isomorphisms and the Five Lemma show that 
$$ 
\cohom_q(G, I_{k+1}\CIc(G)_{\delta}/I_{k}\CIc(G)_{\delta}) \simeq  
I_{k+1} \cohom_q(G, \CIc(G)_{\delta})/I_{k} \cohom_q(G, \CIc(G)_{\delta}).  
$$ 
This is enough to complete the proof. 
\end{proof}

We now study the homology of the subquotients 
$M_k=I_{k+1}\CIc(G)_{\delta}/I_{k}\CIc(G)_{\delta}$ by identifying 
them with induced modules. Let $\Sigma_k$ be a set of representative of
conjugacy classes of standard 
subgroups $S$ such that $F_S \subset V_k\setminus V_{k+1}$ (or, 
equivalently, $F_S^u \subset V_k\setminus V_{k+1}$).

\begin{lemma} \label{Lemma.strata} 
Using the above notation, we have   
$I_k\CIc (G)=\CIc(G\setminus V_k)$ and    
$$ 
I_{k+1}\CIc(G)_{\delta}/I_{k}\CIc(G)_{\delta} \simeq  
\bigoplus\limits_{S \in \Sigma_k} \CIc(F_S^u). 
$$ 
\end{lemma} 
 
\begin{proof} It follows from the definition of $I_k$ that, if   
$f \in I_k\CIc(G)_{\delta}$, then $f$ vanishes in a neighborhood of  
$V_k$. Conversely, if $f$ is in $\CI_c(G \setminus V_k)$, then  
we can find some polynomial $a_{i}$, with $i \le r +k -1$, such  
that $|a_i|$ is bounded from below on the support of $f$  
by, say, $q^{-n}$, then  
$p_nf=f$.  The second isomorphism follows from the first isomorphism 
using \eqref{eq.disc.exact} and Lemma \ref{Lemma.exact}. 
\end{proof}  
 
If $H \subset G$ is a closed  subgroup and $M$ is a smooth (left)  
$H$--module (that is, the stabilizer of each $m \in M$ is an open  
subgroup of $M$), we denote  
\begin{equation}\label{eq.ind} 
	\ind_H^G(M) = \CIc(G) \otimes_H M = \CIc(G) \otimes
	M/\{fh \otimes m - f \otimes hm,\; h \in H\}.
\end{equation} 
Where the right $H$--module structure on $\CIc(G)$ is  
$(fh)(g)=f(gh^{-1})$. Then Shapiro's lemma, see  \cite{Casselman9},
states that 
\begin{equation} 
\label{eq.Shapiro} 
\cohom_k(G, \ind_H^G(M)\otimes \Delta_G) \simeq \cohom_k(H,M). 
\end{equation}
(A proof of Shapiro's Lemma in our setting is contained in the proof of
Theorem \ref{Theorem.Ind}.)

The basic examples of induced modules are obtained from
$H$--spaces. If $X$ is an $H$--space (we agree that $H$ acts on $X$
from the left), then
$$ 
	\CIc((G \times X)/H) \simeq \ind_H^G(\CIc(X) \otimes \Delta_H) 
	\simeq \ind_H^G(\CIc(X)_{\delta})
$$ 
as $G$--modules, where $H$ acts on $G \times X$ by $h(g,x)=(gh^{-1},hx)$. 
For example, Proposition \ref{Prop.stratification} identifies  
$\CIc(F_S^u)$ with an induced module: 
$$
	\CIc(F_S^u) \simeq \ind_{N(S)}^G(\CIc(S^{\rg}\nil S) \otimes
	\Delta_{N(S)})=\ind_{N(S)}^G(\CIc(S^{\rg}\nil S)_{\delta}).
$$ 

Shapiro's lemma is an easy consequence of the Serre--Hochschild spectral
sequence, see \cite{Casselman9}, which states the following. Let $M$ 
be a smooth $G$--module and $H \subset G$ be a normal subgroup. Then
the action of $G$ on $\cohom_q(H,M)$ descends to an action of $G/H$,
and there exists a spectral sequence with $E^2_{p,q}=\cohom_p(G/H,
\cohom_q(H,M))$, convergent to $\cohom_{p+q}(G,M)$.
 
Let $M_k =I_{k+1}\CIc(G)_{\delta}/I_{k}\CIc(G)_{\delta}$, as before.

\begin{proposition} \label{Prop.Gen} 
Using the above notation, we have a natural 
$$ 
\cohom_q(G, M_k) \simeq  
\bigoplus\limits_{\class \in \Sigma_k} 
\CIc(S^{\rg})^{W(S)} \otimes 
\cohom_q(C(S), \CIc(\nil S)_{\delta}) 
$$ 
isomorphism of 
$R^{\infty}$--modules 
\end{proposition}

\begin{proof}  
Let $S$ be a standard subgroup of $G$. Recall first that 
$W(S)=N(S)/C(S)$ is a finite group that acts freely on $S^{\rg}$,
which gives a $N(S)$--equivariant isomorphism 
$$
	\CIc(S^{\rg}\nil S)_{\delta}=\CIc(\nil S)_{\delta} \otimes
	\CIc(S^{\rg}).
$$  
 
Let $M$ be a smooth $N(S)$--module. The Hochschild--Serre spectral 
sequence applied to the module $M$ and the normal subgroup $C(S) 
\subset N(S)$ gives natural isomorphisms 
$$ 
	\cohom_q(N(S),M) \simeq \cohom_0(W(S),\cohom_q(C(S),M)) \simeq 
	\cohom_q(C(S),M)^{W(S)}. 
$$ 
 
Combining these two isomorphisms, we obtain 
\begin{multline*} 
	\cohom_k(G,\CIc(F_S^u)_{\delta}) \simeq
	\cohom_k(G,\ind_{N(S)}^G(\CIc(S^{\rg}\nil S)_{\delta})_{\delta}) \simeq
	\cohom_k(N(S), \CIc(S^{\rg}\nil S)_{\delta})\\ \simeq \big
	(\cohom_k(C(S), \CIc(\nil S)_{\delta}) \otimes \CIc(S^{\rg})
	\big )^{W(S)} \simeq \CIc(S^{\rg})^{W(S)} \otimes
	\cohom_k(C(S), \CIc(\nil S)_{\delta}).
\end{multline*} 
The result then follows from Lemma \ref{Lemma.strata}, which implies
directly that 
$$
	M_k \simeq \oplus_{S \in \Sigma_k} \CIc(F_S^u)_{\delta}.
$$ 
The proof is now complete.
\end{proof} 
 
Combining this proposition with Corollary \ref{Cor.sum}, we obtain the main result 
of this section. Recall that a $p$--adic group $G= \GG(\FF)$ is 
the set of $\FF$--rational points of a linear 
algebraic group $\GG$ defined over a non-archemedean, non-discrete, 
locally compact field $\FF$ of characteristic zero. 
Also, recall that $\nil S$ is the set of
unipotent elements commuting with the standard subgroup $S$, and that
the action of $C(S)$ on $\CIc(\nil S)$ is twisted by the modular function
of $C(S)$, yielding the module $\CIc(\nil S)_{\delta}=
\CIc(\nil S) \otimes \Delta_{C(S)}$.

\begin{theorem} \label{Theorem.Gen} 
Let $G$ be a $p$--adic group. Let  
$\Sigma$ be a set of representative of
conjugacy classes of standard subgroups of 
$S \subset G$ and $W(S)=N(S)/C(S)$, then we have an isomorphism 
$$ 
	\Hd_q(\CIc(G)) \simeq \bigoplus\limits_{S \in \Sigma} 
	\CIc(S^{\rg})^{W(S)} \otimes \cohom_q(C(S), \CIc(\nil S)_{\delta}). 
$$ 
\end{theorem} 
 
{\em Remark.} The isomorphism of the above theorem is not 
natural. A more natural description of $\Hd_q(\CIc(G))$ will 
be obtained in one of the following sections 
by considering higher orbital integrals and  
their Shalika germs

\section{Higher orbital integrals and their Shalika germs}

Proposition \ref{Prop.Gen} of the previous section allows us to
determine the structure of the localized cohomology groups
$\Hd_*(\CIc(G))_{\mfk m}$, where $\mfk m$ is a maximal ideal of
$R^\infty(G)$.  This will lead to an extension of the {\em higher
orbital integrals} introduced by Blanc and Brylinski in
\cite{Blanc-Brylinski}, and to a generalization of some results of
Shalika \cite{Shalika} to higher orbital integrals.  In this way, we
shall also obtain a more natural description of $\Hd_q(\CIc(G))$.
 
First recall the following result.

\begin{proposition} \label{Prop.Rinv} 
Let $G$ be a reductive $p$--adic group over a field of characteristic
$0$, $S\subset G$ be a standard subgroup, and $\gamma \in S^{\rg}$ (that is,
$\gamma$ is a semisimple element such that $C(S)=C(\gamma)$).  
Then there exists a $N(S)$--invariant closed and open neighborhood of 
$\gamma$ in $C(S)$ such that
$$
	G \times U \in (g,h) \to ghg^{-1} \in G
$$ 
defines a homeomorphism of 
$(G \times U)/N(S)$ onto a $G$--invariant, closed and open subset   
$V\subset G$ containing $\gamma$. 
\end{proposition}

\begin{proof} The result follows from Luna's Lemma. For $p$-adic groups,
Luna's Lemma is proved in \cite{Rader-Rallis}, page 109, Properties
``C'' and ``D.'' 
\end{proof}

From this proposition we obtain the following consequences for the 
ring $R^{\infty}(G)$.

\begin{corollary} \label{Cor.Rinv}
Let $U$ and $V$ be as in Proposition \ref{Prop.Rinv} above. 
 
(i) The ring $R^{\infty}(G)$ decomposes as the direct sum $\CI(V)^G
\oplus \CI(V^c)^G$, and $\CI(V)^G \simeq \CI(U)^{N(S)} \subset
\CI(C(S))^{N(S)}=R^{\infty}(C(S))^{W(S)}$. (Here $V^c$ is the
complement of $V$ in $G$.)
 
(ii) For any two semisimple elements $\gamma, \gamma' \in G$, if 
$\phi(\gamma) =\phi(\gamma')$ for all $\phi \in R^{\infty}(G)$, then 
$\gamma$ and $\gamma'$ are conjugated in $G$. 
 
(iii) Let $\mfk m \in R^{\infty}(G)$ be the maximal ideal consisting 
of functions that vanish at a semisimple element $\gamma \in G$. Then
$\mfk p$ is generated by an increasing sequence of projections, and
$M_{\mfk m} \simeq M/{\mfk m}M$, for any $R^{\infty}(G)$--module $M$.   
\end{corollary} 
 
\begin{proof} (i) is an immediate consequence of Proposition 
\ref{Prop.Rinv}. (ii) follows from \cite{Rader-Rallis}, Proposition 2.5.

To prove (iii), observe that the maximal ideal $\mfk m$ is generated
by a sequence of projections $p_n$, that is, $\mfk m = \cup p_n
R^\infty(G)$, with $p_n^2=p_n$.
 
We know from \cite{Rader-Rallis}, Proposition 2.5, that $R^{\infty}(G)$
is isomorphic to $C^{\infty}$, for some locally compact, totally
disconnected topological space $X$. Moreover, if $M$ is a $C^{\infty}(X)$--module
and ${\mfk m}$ is the maximal ideal of functions vanishing at $x_0$, for
some fixed point $x_0 \in X$, then $\CI(X)_{\mfk m} \simeq 
\CI(X)/{\mfk m}\CI(X)$, 
and hence
$$
	M_{\mfk m} = M \otimes_{\CI(X)} \CI(X)_{\mfk m} \simeq 
	M \otimes_{\CI(X)} \CI(X)/{\mfk m}\CI(X) \simeq M/{\mfk m}M.
$$
Since $X$ is metrizable, we can choose a basis $V_n$ of compact
open neighborhoods of $x_0$ in $X$. Then, if we let $p_n$ to be the
characteristic function of $V_n^c$, then $p_n$ are projections
generating ${\mfk m}$. By choosing $V_n$ to be decreasing, we obtain
a decreasing sequence $p_n$. 
\end{proof}

We now consider for each maximal ideal ${\mfk m}\subset  
R^{\infty}=R^{\infty}(G)$ the localization 
$\Hd_q(\CIc(G))_{\mfk m}$.

\begin{proposition} \label{Prop.local} 
Let $\mfk m$ be a maximal ideal of $R^{\infty}(G)$. 
If $\mfk m$ consists of the functions that vanish at the semisimple element 
$\gamma \in G$ and $S\subset G$ is a standard subgroup such that 
$\gamma \in S^{\rg}$, then  
$$ 
	\Hd_q(\CIc(G))_{\mfk m}\simeq \cohom_q(C(S),\CIc(\nil S)_{\delta}). 
$$ 
For all other maximal ideals $\mfk m \subset R^\infty(G)$, we have 
$\Hd_q(\CIc(G))_{\mfk m}=0$. 
\end{proposition}

Note that $\cohom_q(C(S),\CIc(\nil S)) \simeq 
\cohom_q(C(\gamma),\CIc(C(\gamma)_u))$.

\begin{proof}\ The vanishing of $\Hd_*(\CIc(G))_{\mfk m}$ in the second  
half of the proposition follows because $\CI_c(G)_{\mfk m}=0$ in that
case.
 
Assume now that $\mfk m$ consists of functions vanishing at $\gamma$, 
a semisimple element of $G$. The localization functor $V \to 
V_{\mfk m}$ is exact by standard homological algebra.  
The sequence of ideals $(I_k)_{\mfk m}$ is an increasing sequence 
satisfying $(I_0)_{\mfk m}\simeq 0$ and $(I_{m+1})_{\mfk m}\simeq 
\CC$. Choose $k$ such that $(I_k)_{\mfk m} \simeq 0$ and 
$(I_{k+1})_{\mfk m} \simeq \CC$. (This happens if, and only if, $\gamma 
\in V_k \setminus V_{k+1}$.)  It follows that 
$$ 
\cohom_q(G,\CIc(G)_{\delta})_{\mfk m} \simeq 
\cohom_q(G,I_{k+1}\CIc(G)_{\delta}/I_{k}\CIc(G)_{\delta})_{\mfk m}. 
$$ 
 
Since all the isomorphisms of Proposition \ref{Prop.Gen} are compatible  
with this the localization functor, we obtain that  
\begin{multline*} 
	\Hd_q(\CIc(G))_{\mfk m} \simeq
	\cohom_q(G,\CIc(G)_{\delta})_{\mfk m} \\ \simeq
	\bigoplus\limits_{\class \in \Sigma_k} \cohom_q(C(S),
	\CIc(\nil S)_{\delta}) \otimes \big(\CIc(S^{\rg})^{W(S)}/
	{\mfk m}\CIc(S^{\rg})^{W(S)}\big ).
\end{multline*} 
The only quotient $\CIc(S^{\rg})^{W(S)}/ {\mfk m}\CIc(S^{\rg})^{W(S)}$ 
that does not vanish is the one containing (a conjugate of) $\gamma$, 
and then it is isomorphic to $\CC$. This completes the proof. 
\end{proof}

An alternative proof can be obtained by writing  
$$ 
	\cohom_q(G, \CIc(G)_{\delta})_{\mfk m} \simeq \cohom_q(G,
	\CIc(G)_{\delta\mfk m})\simeq \cohom_q(G,
	\CIc(G)_{\delta}/{\mfk m}\CIc(G)_{\delta}),
$$  
and then observing that $\CIc(G)_{\delta}/{\mfk m}\CIc(G)_{\delta} 
\simeq \CIc(\nil S)$, by Corollary
\ref{Cor.Rinv} (iii). However the above proof is more convenient when 
dealing with orbital integrals. See also \cite{N-algg}, first circulated
in 1990 as a preprint of the Mathematical Institute of the Romanian
Academy (INCREST) Nr. 18, March 1990, and where the localization
techniques were first introduced.
 
We now extend the definition of higher orbital integrals introduced 
by Blanc and Brylinski to cover non--regular semisimple elements also.  
Fix a standard subgroup $S \subset G$, and let $k$ be such 
that $S^{\rg} \subset V_k \setminus V_{k+1}$. As in the above proof, 
Proposition \ref{Prop.Gen} gives a natural $R^{\infty}(G)$--linear, 
degree preserving, surjective morphism 
$$ 
	\cohom_*(G,I_{k+1}\CIc(G)_{\delta}/I_k\CIc(G)_{\delta})
	\longrightarrow \CIc(S^{\rg})^{W(S)} \otimes \cohom_*(C(S),
	\CIc(\nil S)_{\delta}),
$$  
and hence a linear map 
\begin{equation}\label{eq.firstS} 
	I_{k+1} \Hd_*(\CIc(G))= \cohom_*(G,I_{k+1}\CIc(G)_{\delta})
	\longrightarrow \CIc(S^{\rg})^{W(S)} \otimes \cohom_*(C(S),
	\CIc(\nil S)_{\delta}).
\end{equation} 

Fix $c \in \cohom^q(C(S), \CIc(\nil S)_{\delta})$ and $\gamma \in S^{\rg}$, 
and let 
$$ 
	\orb_{\gamma,c}=\orb_{\gamma,c}^S: I_{k+1}\Hd_q(\CIc(G))
	\longrightarrow \CC
$$ 
be the evaluation of the map at $\gamma$ and $c$ in \eqref{eq.firstS}. 
We obtain, in particular, that for any $f \in I_{k+1}\Hd_q(\CIc(G))$, 
the function $\gamma \to \orb_{\gamma,c}(f)$ is a locally constant, 
compactly supported function on $S^{\rg}$.  The function $\orb_{\gamma,c}$
can then be 
extended to the whole group $\Hd_q(\CIc(G))$ using a simple 
observation. For any $\gamma \in S^{\rg}$ there exists a locally 
constant function $\phi \in I_{k+1}$ such that $\phi(\gamma)=1$. Then 
let 
$$ 
	\orb_{\gamma,c}(f) := \orb_{\gamma,c}(\phi f),     
$$ 
which is independent of $\phi$. It follows from definition of 
$\orb_{\gamma,c}$ that, for any $f \in \Hd_q(\CIc(G))$, the function 
$\gamma \to \orb_{\gamma,c}(f)$, obtained as above, is a locally 
constant function on $S^{\rg}$, but not necessarily compactly 
supported. We thus obtain the following result.

\begin{proposition}\ Let $S \subset G$ be a standard subgroup. Then 
there exists a degree-preserving, $R^{\infty}(G)$--linear map 
$$ 
	\orb^S:\Hd_*(\CIc(G)) \longrightarrow \CI(S^{\rg})^{W(S)} 
	\otimes \cohom_*(C(S),\CIc(\nil S)_{\delta}), 
$$ 
which is an isomorphism when localized at each maximal ideal $\mfk m 
\subset R^{\infty}(G)$ consisting of functions vanishing at an element 
$\gamma \in S^{\rg}$. 
\end{proposition}

We call the maps $\orb^S$ and $\orb_{\gamma,c}=\orb_{\gamma,c}^S$
``higher orbital integrals'' because they generalize the usual notion
of orbital integral. (If $c$ is a cocycle of dimension $q$, we
call $\orb_{\gamma,c}$ a $q$-higher orbital integral.)
Indeed, assume that $G$ and $C(S)$ are
unimodular, following thorough the identifications in the previous
section, we obtain, for $c_0 = 1\in \cohom^0(C(S),\CIc(\nil S)_{\delta})$ the
evaluation at the identity element $e \in G$, and $f \in \CIc(G)
=\Hd_0(\CIc(G))$, that
$$ 
	\orb_{\gamma,c_0}(f) = \orb_{\gamma,1}(f) = 
	\int_{G/C(S)}f(g\gamma g^{-1})d\overline{g}, 
$$ 
where $d\overline{g}$ is the induced measure on $G/C(S)$. 
 
If $\gamma \in G$ is a semisimple element and $S$ is a standard 
subgroup of $G$ such that $C(\gamma) = C(S)$, (\ie $\gamma \in  
S^{\rg}$), then restriction at $\gamma$ defines a map  
$$ 
	\orb_\gamma=\orb_\gamma^S: \Hd_*(\CIc(G)) \to
	\cohom_*(C(S),\CIc(\nil S)_{\delta})
$$ 
such that $c(\orb_\gamma(f)) = \orb_{\gamma,c}(f)$, for all 
$c \in \cohom^q(C(S),\CIc(\nil S)_{\delta})$. 
 
A word on notation, whenever we write $\orb_{\gamma,c}^S$ or
$\orb_{\gamma}^S$, we assume that $\gamma \in S^{\rg}$, which actually
determines $S$. This means that we can omit $S$ from notation.
However, if we want to write that $\orb_{\gamma,c}=\orb_{\gamma,c}^S$
is obtained by specializing
$$
	\orb^S:\Hd_*(\CIc(G)) \to \CI(S^{\rg})^{W(S)} \otimes
	\cohom_*(C(S),\CIc(\nil S)_{\delta})
$$  
to a point $\gamma \in S^{\rg}$ and then by evaluating at $c$, that is that 
$$ 
	\orb_{\gamma,c}^S(f)=\langle \orb^S(f)(\gamma), c \rangle, 
$$ 
then it is obviously more convenient to include $S$ in the 
notation. 
 
Let $\gamma \in G$ be a semisimple element. We want now to investigate
the behavior  %behavior behavior behavior behavior behavior 
of the orbital integrals $\orb_{g,c}$ with $g$ in a
small neighborhood of $\gamma$. Fix a standard subgroup $S \subset G$
such that $\gamma$ is in the closure of $Ad_{G}(S^{\rg})$, but is not
in $Ad_G(S^{\rg})$, and a class $c \in \cohom^q(C(S),\CIc(\nil
S))$. More precisely, we want to study the germ of the function $g \to
\orb_{g,c}(f)$ at an element $\gamma$, where $f \in \Hd_q(\CIc(G))$ is
arbitrary. The germ of a function $h$ at $\gamma$ will be denoted
$h_\gamma$.
 
The following theorem extends one of the basic properties of 
Shalika germs from usual orbital integrals to higher orbital 
integrals.

\begin{theorem} \label{Theorem.Shalika} 
Let $S \in G$ be a standard subgroup $\gamma \in S$ an element in the
closure of $S^{\rg}$, but $\gamma \not \in S^{\rg}$.  Then there
exists a degree preserving linear map
$$ 
	\sigma_{\gamma}^S: \cohom_*(C(\gamma),\CIc(C(\gamma)_u)_{\delta}) 
	\rightarrow \CI(S^{\rg})^{W(S)}_{\gamma}  
	\otimes \cohom_*(C(S),\CIc(\nil S)_{\delta}),  
$$ 
such that  
$$ 
\orb^S(f)_{\gamma}=\sigma^S_\gamma(\orb_{\gamma}(f)),
$$ 
for all $f \in \Hd_*(\CIc(G))$. 
\end{theorem}

Note that, in the notation for the maps $\sigma_{\gamma}^S$, the standard 
subgroup $S$ is no longer determined by $\gamma$.

\begin{proof} By the definition of the localization of 
a module, the map 
$$
	\orb^S :\Hd_*(\CIc(G)) \to \CI(S^{\rg})^{W(S)}_{\gamma}\otimes
	\cohom_*(C(S),\CIc(\nil S)_{\delta})
$$  
factors through a map  
$$ 
	F: \Hd_*(\CIc(G))_{\gamma} \to \CI(S^{\rg})^{W(S)}_{\gamma}
	\otimes \cohom_*(C(S),\CIc(\nil S)_{\delta}).
$$  
Since $\orb_\gamma : \Hd_*(\CIc(G))_{\gamma} \to  
\cohom_*(C(\gamma), \CIc(C(\gamma)_u)_{\delta})$ is an isomorphism by 
Proposition \ref{Prop.local}, we may define  
$$ 
	\sigma_{\gamma}^S = F \circ \orb_{\gamma}^{-1}, 
$$ 
and all desired properties for $\sigma_{\gamma}^S$ will be satisfied. 
\end{proof}

Let $\gamma \in S \setminus S^{\rg}$ be such that $\gamma$ is in the
closure of $S^{\rg}$, as above, and also let $c \in
\cohom^q(C(S),\CIc(\nil S)_{\delta})$. Then a consequence of the above
theorem, Theorem \ref{Theorem.Shalika}, is that the germ at $\gamma$
of the higher orbital integrals $\orb_{g,c}^S$ depends only on
$\orb_\gamma$. More precisely, if $g \in S^{\rg}$, $f \in
\Hd_q(\CIc(G))$, and we regard $\orb_{g,c}^S(f)$ as a function of $g$,
then its germ at $\gamma$, denoted $\orb_{g,c}^S(f)_{\gamma}$, is
given by
\begin{equation} 
	\orb_{g,c}^S(f)_{\gamma} = \langle
	\sigma_{\gamma}^S(\orb_{\gamma}(f)), c \rangle.
\end{equation}  
 
This observation allows us to relate Theorem \ref{Theorem.Shalika}
with results of Shalika \cite{Shalika} and Vigneras
\cite{Vigneras}. So assume now that $G$ is reductive and let $\xi_i
\in \cohom_0(C(\gamma),\CIc(C(\gamma)_u)_{\delta})$ be the basis dual
to the basis of $\cohom^0(C(\gamma),\CIc(C(\gamma)_u)_{\delta})$ given by the
orbital integrals over the orbits of $\gamma u$, for $u$ nilpotent in
$C(\gamma)$. If we let $F_i = \sigma_{\gamma}(\xi_i)$, then we recover
the usual definition of Shalika germs. Due to this fact, we shall call
the maps $\sigma_{\gamma}^S$ introduced in Theorem
\ref{Theorem.Shalika} the {\em higher} Shalika germs.
 
We can now characterize the range of the higher orbital integrals. 
Combining all higher orbital integrals for $S \subset G$ ranging 
through a set $\Sigma$ of representatives of standard subgroups of $G$,  
we obtain a map 
$$ 
\orb: \Hd_*(\CIc(G)) \to \bigoplus\limits_{S \in \Sigma} 
\CI(S^{\rg})^{W(S)} \otimes \cohom_*(C(S),\CIc(\nil S)). 
$$

\begin{theorem}\ Let $\Sigma$ be a set of representatives of 
standard subgroups of $G$ and $\sigma_{\gamma}^S$ be the maps 
introduced in Theorem \ref{Theorem.Shalika} for $\gamma \in 
\overline{S^{\rg}} \setminus S^{\rg}$. Also, let  
$$ 
	{\mathcal F} \subset \bigoplus\limits_{S \in \Sigma}  
	\CI(S^{\rg}) \otimes \cohom_*(C(S),\CIc(\nil S)_{\delta}), 
$$  
be the space of sections $\xi$ satisfying $\xi_{\gamma}= 
\sigma_{\gamma}^S(\xi(\gamma))$. 
Then $\orb$ establishes a $R^{\infty}(G)$--linear isomorphism  
$$ 
	\orb:\Hd_*(\CIc(G)) \longrightarrow {\mathcal F}. 
$$ 
\end{theorem}

\begin{proof} Note first that the map $\orb$ is well defined,
that is, that its range is contained in ${\mathcal F}$, 
by Theorem \ref{Theorem.Shalika}.  
 
To prove that $\orb$ is an isomorphism, filter both  
$\Hd_*(\CIc(G))$ and ${\mathcal F}$ by the subgroups $I_k\Hd_*(\CIc(G))$ 
and, respectively, by $I_k{\mathcal F}$,  
using the ideals $I_k$ introduced in Section \ref{Sec.Homology}. 
Since $\orb$ is $R^{\infty}(G)$--linear, it preserves this filtration 
and induces maps 
$$
I_{k+1}\Hd_*(\CIc(G))/I_k\Hd_*(\CIc(G)) 
\to I_{k+1}{\mathcal F}/I_k{\mathcal F}.
$$
These maps are, by construction, exactly the isomorphisms of Proposition 
\ref{Prop.Gen}. Standard homological algebra then implies that  
$\orb$ itself is an isomorphism, as desired. 
\end{proof}

A consequence of the above result the following ``density'' 
corollary.

\begin{corollary}\ Let $a \in \Hd_q(\CIc(G))$. If all $q$-higher
orbital integrals of $a$ vanish, then $a=0$.
\end{corollary}

We shall also need certain specific cocycles below. Let $\tau_0$ be the
trace $\tau_0(f)=f(e)$ on $\CIc(G)$, $G$ unimodular, obtained by 
evaluating $f$ at the identity $e$ of $G$. Let $G_0$ be the kernel of
all unramified characters of $G$. Then $G/G_0 \simeq \ZZ^r$, where
$r$ is the rank of a split component of $G$. Let $p_j : G \to \ZZ$
be the morphisms obtained by considering the $j$th component of $\ZZ^r$.
Then
$$
	\delta_j(f)(g) = p_j(g) f(g)
$$
defines a derivation of $\CIc(G)$. Moreover, we can identify
$\cohom^*(G)$ with $\Lambda^*\CC^r$,
the exterior algebra with generators $\delta_1, \ldots, \delta_r$.
Fix $c \in \cohom^*(G)$. We can assume that $c = \delta_1 \wedge
\ldots \wedge \delta_q$, and then we define the map
$D_c : \CIc(G)^{\otimes q+1} \to \CIc(G)$ by the formula
\begin{equation}\label{eq.Dc}
	D_c(f_0,\ldots,f_q) = (q!)^{-1}\sum_{\sigma \in S_q}
	\epsilon(\sigma) f_0 \delta_{\sigma(1)}(f_1) 
	\delta_{\sigma(2)}(f_2) \ldots \delta_{\sigma(q)}(f_q) 
\end{equation}
Then $\tau_c=\tau_0\circ D_c(f_0,\ldots, f_q)$
defines a Hochschild $q$ cocycle on $\CIc(G)$, and 
\begin{equation}\label{eq.cc}
	\tau_c = \orb_{e,c},
\end{equation}
if we naturally identify $c$ with an element of the cohomology
group $\cohom^*(G, \CIc(G_u))$.

\section{The cohomology of the unipotent variety}

It follows from the main result of the first section, Theorem
\ref{Theorem.Gen}, that, in order to obtain a more precise description of
the Hochschild homology of $\CIc(G)$, we need to understand the
continuous cohomology of the $H$--module $\CIc(H_u)_{\delta}$, where
$H$ ranges through the set of centralizers of standard subgroups of $G$ and
$H_u$ is the variety of unipotent elements in $H$. (We call the
variety $H_u$ the {\em unipotent variety} of $H$.)  Since the
cohomology groups $\cohom_k(H, \CIc(H_u)_{\delta})$ depend only on
$H$, it is enough to consider the case $H=G$.  In this section we gather
some results on the groups $\cohom_q(G, \CIc(G_u)_{\delta})$.
 
We first need to recall the computation of the groups
$\cohom_*(G)=\cohom_*(G,\CC)$, \cite{Borel-Wallach}. 
More generally, we also need to compute
$\cohom_*(G,\CC_{\chi})$, where $\chi:G \to \CC^*$ is a character of
$G$ and $\CC_{\chi}=\CC$ as a vector space, but with $G$--action given
by the character $\chi$.
 
Assume first that $G=S$ is a commutative $p$--adic group, and let
$S_0$ be the union of all compact-open subgroups of $S$. Then $S_0$ is
a subgroup of $S$ and $S/S_0$ is a free abelian subgroup, whose rank
we denote by $\rk(S)$. For this group we then have
$$ 
\cohom_k(S) \simeq \cohom_k(S/S_0,\CC)   
\simeq \Lambda^k \CC^{\rk(S)}.  
$$ 
 
For an arbitrary $p$--adic group $G$, we may identify the cohomology
groups $\cohom_q(G)$ with those of a commutative $p$--adic group.
Indeed, if $G^0$ is the connected component of $G$ (in the sense of
algebraic groups) then $G/G^0$ is finite, and hence $\cohom_q(G)
\simeq \cohom_q(G^0)$, by the Hochschild--Serre spectral
sequence. This tells us that we may assume $G$ to be connected as an
algebraic group. Choose then a Levi decomposition $G=MN$, where $N$ is
the unipotent radical of $G$, $M$ is a reductive subgroup, uniquely
determined up to conjugation, and the product $MN$ is a semidirect
product.  Since $\cohom_q(N)=0$ for $q >0$, it follows that
$\cohom_q(G) \simeq \cohom_q(M)$. Let $M_1 \subset M$ be the
commutator subgroup of $M$, which is also a $p$--adic group, see
\cite{Borel1}. The cohomology groups $\cohom_q(M_1)$ were computed in
\cite{Blanc-Brylinski}, Proposition 6.1, page 316, and
\cite{Borel-Wallach} and they also vanish for $q>0$ (because the
fundamental domain of the building of $G$ is a simplex).  All in all,
we obtain that
$$\cohom_q(G) \simeq \cohom_q(M) \simeq \cohom_q 
(M^{\operatorname{ab}}), 
$$ 
where $M^{\operatorname{ab}}=M/M_1$ is the abelianization of $M$. 
 
We summarize the above discussion in the following well known  
statement. 
 
\begin{lemma} Let $G$ be a $p$--adic group, not necessarily reductive, 
and let $r$ be the rank of a split component of the reductive  
quotient of $G$. Then 
$$ 
\cohom_q(G) =\cohom_q(G,\CC)\simeq \Lambda^q \CC^r, 
$$ 
and $\cohom_q(G,\CC_{\chi})=0$, if $\chi$ is 
a nontrivial character of $G$.
\end{lemma}

We continue with a few elementary remarks on $\cunip {k}{G}$.   
% 
%     $\cunip {k}{G} = 
% 
%     $\cohom_k(G,\CI_c(G_u)_{\delta})$ 
% 

\begin{remark} If $G_1 \to G$ is a surjective morphism with finite 
kernel $F$, then there exists a natural homeomorphism $G_{1u}\simeq G_u$ 
of the unipotent varieties of the two groups. Since the kernel $F$ 
acts trivially on $G_{1u}$, using the Hochschild-Serre spectral 
sequence we obtain an isomorphism 
\begin{equation} 
\cunip {k}{G_1} \simeq \cohom_k(G_1, \CI_c(G_u)) \simeq 
\cunip {k}{G}. 
\end{equation} 
\end{remark} 

\begin{remark} If $G \subset G_1$ is a normal $p$--adic subgroup with 
$F \simeq G_1/G$ finite, then we again have a natural homeomorphism 
$G_{1u} = G_u$. This gives 
\begin{equation} 
	\cunip {k}{G_{1}} \simeq \cohom_k(G_1, \CI_c(G_u)_{\delta})
	\simeq {\cunip {k}{G}}^{F},
\end{equation} 
using once again the Hochschild-Serre spectral sequence. In 
particular, if the characteristic morphism $F \to 
\operatorname{Aut}(G)/ \operatorname{Inn}(G)$ is trivial, then we get 
a natural isomorphism $\cunip {k}{G} \simeq \cunip {k}{G_{1}}$. 
\end{remark} 

\begin{remark} If $G = G' \times G''$, then $G_u = G_u' \times G_u''$ 
naturally, and hence $\CIc(G_u) \simeq \CIc(G_u') \otimes
\CIc(G_u'')$.  This gives
\begin{equation} 
	\cunip {k}{G} \simeq \bigoplus\limits_{i+j = k} 
	\cohom_i(G', \CIc(G_u')) \otimes  \cohom_j(G'', \CIc(G_u''))   .
\end{equation} 
\end{remark} 

\begin{remark} If $Z$ is a commutative $p$--adic group of split rank 
$r$, then $$\cunip {k}{Z} \simeq \CIc(Z_u) \otimes \Lambda^k \CC^r.$$ 
\end{remark} 

\begin{remark} The above isomorphisms reduce the computation of 
$\cunip {k}{G}$ for $G$ reductive, to the computation of the cohomology 
groups corresponding to its semisimple quotient $H:=G/Z(G)$: 
\begin{equation} 
	\cohom_k(G, \CIc(G_u)_{\delta}) = \cohom_k(G, \CIc(G_u))\simeq
	\bigoplus\limits_{i + j = k} \cunip {i}{H} \otimes \Lambda^j
	\CC^r,
\end{equation} 
where $r$ is the rank of a split component of $G$. 

Let $\tau_0$ be the
trace obtained by evaluating at the identity. Using $\tau_0$, we obtain
an injection $H^j(G) \ni c \to \tau_0 \otimes c \in H^j(G,\CIc(G_u))$.
\end{remark} 
 
In order to obtain more precise results on $\cohom_*(G,
\CIc(G_u)_{\delta})$, we need to take a closer
look at the structure of $\CIc(G_u)$ as a $G$--module.  For a
$G$--space $X$, we denote by $\langle X \rangle$ the quotient space
$X/G$ with the induced topology, which may be non-Hausdorff.  Thus
$\langle G_u \rangle$ is the set of unipotent conjugacy classes of
$G$.
 
Assume now that $\langle G_u \rangle$ is a finite set. (This happens
for example if $G$ is reductive, because the ground field $\FF$ has
characteristic zero.) Then the space $G_u$ can be written as an
increasing union of open $G$--invariant sets $U_l \subset G_u$,
$U_{-1} = \emptyset$, such that each difference set $U_l \setminus
U_{l-1}$ is a disjoint union of open and closed $G$--orbits,
\begin{equation}\label{eq.orbits} 
U_l \setminus U_{l-1} = \cup X_{l,j}.
\end{equation} 
A filtration $U_l$ with these properties will be called ``nice.'' 
There may be several nice filtrations of $G_u$. 
 
A nice filtration of $G_u$, as above, gives rise, by standard
arguments, to a spectral sequence converging to $\cunip {k}{G}$, as
follows. First, let $\langle g \rangle \in \langle G_u \rangle$ be the
orbit through an element $g \in G_u$. Also, let $C(g)$ denote the
centralizer of $g \in G_u$ and $r_g$ denote the rank of a split
component of $C(g)$ if $C(g)$ is unimodular, $r_g=0$ otherwise.  This
definition of $r_g$ is justified by $\cohom_k (C(g), \Delta_{C(g)})
\simeq \Lambda^k \CC^{r_g}.$

\begin{proposition}\label{Prop.Spectral} \
Let $G$ be a  $p$--adic group with finitely many  
unipotent orbits (\ie $\langle G_u \rangle$ is finite). Then, for 
any nice filtration $(U_l)$ of $G_u$ by open $G$--invariant subsets, 
there exists a natural spectral sequence with  
$$ 
	E^2_{p,q} = \bigoplus_{\langle u \rangle \in \langle
	U_{p}\setminus U_{p-1}\rangle} \Lambda^{p+q}\CC^{r_u},
$$ 
convergent to $\cunip {p+q}{G}$. 
\end{proposition}

\begin{proof}\ The argument is standard and goes as follows. Recall first 
that any filtration $0 = F_0 \subset F_1 \subset \ldots \subset F_N 
=\CIc(G_u)_{\delta}$ by $G$--submodules gives rise to a spectral sequence with 
$E^1_{p,q}=\cohom_{p+q} (G, F_{p}/F_{p-1})$, convergent to $\cunip 
{p+q}{G}$. 
 
Now, associated to the open sets $U_l$ of a nice filtration, there 
exists an increasing 
filtration $F_l=\CIc(U_l)_{\delta} \subset \CIc(G_u)_{\delta}$ by 
$G$--submodules such that 
$$ 
\CIc(U_l)_{\delta}/\CIc(U_{l-1})_{\delta}   \simeq \bigoplus\limits_{j} 
\CIc(X_{l,j})_{\delta}, 
$$ 
where each $X_{l,j}$ is the orbit of a unipotent element, and $U_l
\setminus U_{l-1}$ has the topology given by the disjoint union of the
orbits $X_{l,j}$.  Fix $l$ and $j$, and let $u$ be a unipotent element
in $X_{l,j}$ (so that then $X_{l,j}$ is the orbit through $u$), which
implies that $\CIc(X_{l,j})\simeq
\ind_{C(u)}^{G}(\Delta_{C(u)})$. Finally, from Shapiro's lemma we
obtain that
$$ 
	\cohom_k(G,\CIc(X_{l,j})_{\delta}) \simeq
	\cohom_k(C(u),\Delta_{C(u)}) \simeq \Lambda^k \CC^{r_u},
$$ 
and this completes the proof.
\end{proof}

We expect this spectral sequence to converge for $G$ reductive.  This
is the case, for example for $G=GL_n(\FF)$ and for $SL_n(\FF)$. See
Section \ref{Sec.Examples}.  The convergence of the spectral sequence
implies, in particular, the convergence of the orbital integrals of
unipotent elements in reductive groups (which is a well known fact due
to Deligne and Rao \cite{Rao}).  In general, the convergence of the
spectral sequence of the above proposition can be interpreted to
represent the convergence of ``higher orbital integrals.''

\section{Induction and the unipotent variety}

We assume from now on that $G$ is reductive, and we fix a parabolic
subgroup $P \subset G$, $P \not = G$, and a Levi subgroup $M \subset
P$, so that $P=MN$, where $N$ is the unipotent radical of $P$, and the
product is a semidirect product. In this section, we relate $\cunip
{*}{G}$ to the groups, $\cunip *{P}$ and $\cunip *{M}$.  He have
considered non-unimodular subgroups in the previous sections in order
to be able to handle subgroups like $P$.
 
Let $K$ be a ``good'' maximal compact subgroup of $G$ (see \cite{HC},
Theorem 5), so  that $G = KP$. This decomposition shows that the map 
$$ 
	K \times P \ni (k,p) \to kpk^{-1} \in G 
$$  
is proper, and hence the map $G \times_P P:=(G \times P)/P \ni (g,p) 
\to gpg^{-1} \in G$ is also proper. This gives a map 
$$
	\CIc(G)_{\delta}=\CIc(G) \to \CIc(G \times_P P)
	\simeq \ind_P^G(\CIc(P)
	\otimes \Delta_P)=\ind_P^G(\CIc(P)_{\delta})
$$ 
of $G$-modules. This map of $G$--modules and the standard identification 
of Hochschild homology with continuous cohomology, equation 
\eqref{eq.BlBr}, then give a morphism 
\begin{equation}\label{eq.ind.morphism} 
	\ind_P^G :\Hd_*(\CIc(G))\longrightarrow \Hd_*(\CIc(P)), 
\end{equation} 
defined as the composition of the following sequence of morphisms 
\begin{multline*} 
	\Hd_*(\CIc(G)) \simeq \cohom_*(G,\CIc(G)_{\delta}) \longrightarrow 
	\cohom_*(G, \ind_P^G(\CIc(P)_{\delta})\otimes \Delta_G) \\ 
	\simeq \cohom_k(P,\CIc(P)_{\delta}) \simeq \Hd_*(\CIc(P)) 
\end{multline*} 
of Hochschild homology groups.  The main result of this section 
states that $\ind_P^G$ is induced by a morphism of algebras, which we 
now proceed to define. 
 
Let $dk$ be the normalized Haar measure on the maximal compact 
subgroup $K$, normalized such that $K$ has volume $1$.  The 
composition of kernels 
$$ 
	T_1T_2(k_1,k_2) = \int_{G/P} T_1(k_1,k) T_2(k,k_2) dk 
$$ 
defines on $\CI(K \times K)$ an algebra structure. Let
\begin{equation}\label{eq.morphism} 
\phi_P^G: \CIc(G) \longrightarrow \CI(K \times K) \otimes \CIc(P), 
\end{equation} 
be defined by $\phi_P^G(f)(k_1,k_2,p)=f(k_1 p k_2^{-1})$. 
 
Recall \cite{HC} that the push-forward of the product $dpdk$ of Haar
measure on $P \times K$, via the multiplication map $P \times K \ni
(p,k) \to pk \in G$, is a left invariant measure on $G$, and hence a
multiple $\lambda dg$ of the Haar measure $dg$ on $G$. Suppose that
the measure $dk$ of $K$ is the restriction of $dg$ to $K$, and has
total mass 1.  Then the Haar measures on $G$ and $P$ will be called
{\em compatible} if $\lambda = 1$. We shall need the following result of Harish
Chandra (implicitly stated in \cite{Silberger}):

\begin{lemma}\ Suppose the Haar measures on $G$ and $P$ are compatible.  
Then the linear map $\phi_P^G$, defined above in equation
\eqref{eq.morphism}, is a morphism of algebras.
\end{lemma}

\begin{proof}\ The product on $\CI(K \times K) \otimes \CIc(P)= 
\CIc(K \times K \times P)$ is given by the formula
$$ 
	(h_1 h_2)(k_1,k_2,p)= \int_K \int_P h_1(k_1,k,q)
	h_2(k,k_2,q^{-1}p)dq dk.
$$ 
Let $*$ denote the multiplication (\ie convolution product) on 
$\CIc(G)$. Thus, we need to prove that 
$$ 
	f_1 *f_2(k_1pk_2^{-1})= \int_K \int_P f_1(k_1qk^{-1})
	f_2(kq^{-1}pk_2^{-1})dq dk,
$$ 
for all $f_1,f_2 \in \CIc(G)$. Consider the map $P \times K \ni (q,k) 
\to g:=qk^{-1} \in G$, and let $d\mu$ be the push-forward of the measure 
$dq dk$. Then the right-hand side of the above formula becomes 
$$ 
	\int_K \int_P f_1(k_1qk^{-1}) f_2(kq^{-1}pk_2) dq dk = \int_G f_1(k_1g) 
	f_2(g^{-1} p k_2^{-1})d \mu(g). 
$$ 
We know that $d\mu = dg$, by assumptions (see the 
discussion before the statement of this lemma), and then 
$$ 
\int_G f_1(k_1g) 
f_2(g^{-1} p k_2^{-1})d \mu(g)= 
\int_G  f_1(g) f_2(g^{-1}k_1 p k_2^{-1}) d\mu(g)= 
f_1* f_2(k_1 p k_2^{-1}),
$$ 
by the invariance of the Haar measure. The lemma is proved. 
\end{proof} 
 
The trace $\CI(K \times K) \to \CC$ induces an isomorphism  
$$
\tau:\Hd_*(\CI(K \times K) \otimes \CIc(P))\simeq \Hd_*(\CIc(P)).
$$  
Explicitly, this isomorphism is given at the level of chains by 
\begin{multline*} 
\tau(f_0\otimes f_1 \otimes \ldots \otimes f_q)(p_0,p_1, 
\ldots,p_q) \\ 
:= \int_{K^{q+1}} f_0(k_0,k_1,p_0)f_1(k_1,k_2,p_1)  
\ldots f_q(k_q,k_0,p_q) dk_0 \ldots dk_q. 
\end{multline*} 
This isomorphism combines with $\phi_P^G$ to give a morphism 
\begin{equation}\label{eq.above} 
	(\phi_P^G)_*: \Hd_*(\CIc(G)) \longrightarrow \Hd_*(\CIc(P)). 
\end{equation}

\begin{theorem} \label{Theorem.Ind}
Let $P$ be a parabolic subgroup of a reductive $p$--adic group 
$G$. Consider the morphisms $(\phi_P^G)_*$ and 
$\ind_P^G:\Hd_*(\CIc(G)) \longrightarrow \Hd_*(\CIc(P))$, defined 
above (Equations \eqref{eq.above} and \eqref{eq.ind.morphism}). 
Then $(\phi_P^G)_*=\ind_P^G$. 
\end{theorem}

\begin{proof} Let $M_1$ and $M_2$ be two left $G$--modules. We 
can regard $M_1$ as a right module, and then the tensor product $M_1 
\otimes_G M_2$ is the quotient of $M_1 \otimes M_2$ by the group 
generated by the elements $gm_1 \otimes gm_2 - m_1 \otimes 
m_2$, as before. Alternatively, we can think of $M_1 \otimes_G M_2$ as $(M_1 
\otimes M_2) \otimes_G \CC$.  This justifies the notation $f\otimes_G 
1$ for a morphism $M_1 \otimes_G M_2 \to M_1' \otimes_G M_2'$ induced 
by a morphism 
$$
	f =f_1 \otimes f_2 : M_1 \otimes M_2 \to M_1' \otimes M_2'.
$$ 
  
We shall prove the theorem by an explicit computation. To this end, we
shall use the results and notation ($h_G$ and $\tilde h_G = h_G
\otimes_G 1$) of Lemma \ref{Lemma.BlBr}.

By direct computation, we see that the morphism  
$$
	\tau \circ \phi_P^G :\CIc(G)^{\otimes q+1} \to
	\CIc(P)^{\otimes q+1}
$$  
between Hochschild complexes, is given by the formula
\begin{multline} \label{eq.phiPG} 
	\tau \circ \phi_P^G(f) (p_0,p_1,\ldots,p_q)=\\ \int_{K^{q+1}}
	f(k_0p_0k_1^{-1},k_1p_1k_2^{-1},\ldots,k_qp_qk_{0}^{-1})
	dk_0dk_1 \ldots dk_q.
\end{multline} 
 
We now want to realize the map $\ind_P^G:\Hd_*(\CIc(G))
\longrightarrow \Hd_*(\CIc(P))$, at the level of complexes. In the
process, it will be convenient to identify the smooth $G$-module
$\CIc((G \times P)/P) \simeq \ind_P^G(\CIc(P)_{\delta})$ 
with a subspace of the space of functions on $G \times P$, using
the projection $G \times P \to (G \times P)/P$.
 
Consider the $G$--morphism 
$$	
	l:{\mathcal B}_q(G) \otimes \CIc(G) \to {\mathcal B}_q(G)
	\otimes \ind_P^G(\CIc(P)_{\delta})
$$ 
induced by the morphism 
$$
	\CIc(G) \to \ind_P^G(\CIc(P)_{\delta})\subset \CI(G \times P).
$$
Explicitly, 
$$
	l(f)(g_0,g_1,\ldots,g_q,g,p)= f(g_0,g_1,\ldots,g_q,gpg^{-1}).
$$
Then the resulting morphism
$$
	l \otimes_G 1: \cohom_q(G,\CIc(G)_{\delta}) 
	= \cohom_q(G,\CIc(G)_{\delta}) \to 
	\cohom_q(G,\ind_P^G(\CIc(P)_{\delta}))
$$
is the morphism $\cohom_q(G,\CIc(G)_{\delta}) \to
\cohom_q(G,\ind_P^G(\CIc(P)_{\delta}))$ on homology corresponding
to the $G$--morphism $\CIc(G) \to \ind_P^G(\CIc(P)_{\delta})$.
 
The $G$--morphism
\begin{multline} \label{eq.a}
	r:{\mathcal B}_q(G) \otimes \ind_P^G(\CIc(P)_{\delta}) \to
	\ind_P^G({\mathcal B}_q(P) \otimes \CIc(P)_{\delta}) \\ = \CIc(G)
	\otimes_P ({\mathcal B}_q(P) \otimes \CIc(P)_{\delta}) \subset \CIc(G
	\times P^{q+2}),
\end{multline} 
given by the formula 
$$ 
	r(f)(g,p_0,p_1,\ldots, p_q,p)= 
	\int_{K^{q+1}}f(gp_0k_1^{-1},gp_1k_2^{-1}, 
	\ldots,gp_nk_{0}^{-1},g,p)dk, 
$$ 
($dk = dk_0 \ldots dk_q$)
is well defined and induces an isomorphism in homology, because the 
only nonzero homology groups are in dimension $0$, and they are both 
isomorphic to $\ind_P^G(\CIc(P)_{\delta})$. We have an isomorphism  
$$ 
\chi: \ind_P^G({\mathcal B}_q(P) \otimes \CIc(P)_{\delta}) \otimes_G \CC
\to ({\mathcal B}_q(P) \otimes \CIc(P)_{\delta}) \otimes_P \CC, 
$$ 
of complexes. This shows that the homology of  
the second complex in \eqref{eq.a}
is isomorphic to $\cohom_q(P, \CIc(P)_{\delta})$, and that the map 
induced on homology, that is 
$$
	\chi (r \otimes_G 1) : \cohom_q(G,\ind_P^G(\CIc(P)_{\delta})) \to
	\cohom_q(P,\CIc(P)_{\delta}),
$$ 
is the Shapiro isomorphism.
 
Recall now that the isomorphism 
$\cohom_q(G,\CIc(G)_{\delta}) \simeq \Hd_q(\CIc(G))$ is
induced by the morphism of complexes $\tilde{h}_q$ defined
in Lemma \ref{Lemma.BlBr}, Equation \eqref{eq.BlBr}.
From the definition of the morphism $\ind_P^G : \Hd_*(\CIc(G))
\to \Hd_*(\CIc(P))$ and the above discussion, 
we obtain the equality of the morphisms $\cohom_q(G,\CIc(G)_{\delta}) \to 
\cohom_q(P,\CIc(P)_{\delta})$ induced by $\chi(rl \otimes_G 1)$  
and $\tilde{h}_P^{-1}\circ \ind_P^G \circ \tilde{h}_G$. Thus, in order 
to complete the proof, it would be enough to check that 
$\tilde{h}_P\circ \chi(rl \otimes_G 1) = \tau \circ \phi_P^G 
\circ \tilde{h}_G$ at the level of complexes. Let 
$$
	t : {\mathcal B}_q(G) \otimes \ind_P^G(\CIc(P)_{\delta}) 
	\to {\mathcal B}_q(G) \otimes_G \ind_P^G(\CIc(P)_{\delta})
$$ 
be the projection. Since $h_G$ is surjective, 
it is also enough to check that  
$\tilde{h}_P\chi (r \otimes_G 1) t l=\tau \phi_P^G h_G$.

Let  
$$ 
	r'(f)(p_0,p_1,\ldots, p_q,p)= 
	\int_{K \times K^{q+1}}f(\,k'p_0k_1^{-1},\,k'p_1k_2^{-1},\, 
	\ldots,\,k'p_nk_{0}^{-1},\,k',\,p)dk'dk, 
$$ 
where $dk = dk_0 \ldots dk_q$, as before.  Then $r'$ induces a
morphism 
$$
	r': {\mathcal B}_q(G) \otimes \ind_P^G(\CIc(P)_{\delta}) \to
	{\mathcal B}_{q}(P) \otimes \CIc(P)_{\delta}
$$ 
of complexes satisfying $h_P \circ r'= \tilde{h}_P \chi(r \otimes_G 1)
t$. Directly from the definitions we obtain then that $h_P\circ r'
\circ l = \tau \circ \phi_P^G \circ h_G$.  This completes the proof.
\end{proof}

For simplicity, we have stated and proved the above result only
for $G$ reductive, however, it extends to general $G$ and $P$ such
that $G/P$ is compact, by including
the modular function of $G$, where appropriate. 
 
In order to better understand the effect of the morphism  
$$\ind_P^G=(\phi_P^G)_*:\Hd_*(\CIc(G)) \longrightarrow \Hd_*(\CIc(P)),$$ 
it is sometimes useful to look at its action on the geometric fibers 
of the group $\Hd_*(\CIc(G))$. This is especially useful because 
the action on the geometric fibers also recovers the 
classical results on the characters of induced representations. 
 
First we observe that restriction defines a morphism  
$\rho_P^G : R^{\infty}(G) \to R^{\infty}(P)$. In case the group $G$ is  
reductive and $M$ is 
a Levi component of the parabolic subgroup $P$, we also have 
$R^{\infty}(P) \simeq R^{\infty}(M)$.

\begin{lemma}\ Let $P$ be a parabolic subgroup of a reductive  
$p$--adic group $G$, and let $\rho_P^G : R^{\infty}(G) \to R^{\infty}(P)$ be the 
morphism induced by restriction, used to define a  
$R^{\infty}(G)$--module structure on $\Hd_*(\CIc(P))$. Then 
$$
	\ind_P^G : \Hd_*(\CIc(G)) \longrightarrow \Hd_*(\CIc(P))
$$ 
is $R^{\infty}(G)$--linear, 
in the sense that $\ind_P^G(f \xi) = \rho_P^G(f) \ind_P^G(\xi)$, for 
all $f \in R^{\infty}(G)$ and all $\xi \in \Hd_*(\CIc(G))$. 
\end{lemma}

\begin{proof}\ The result of the lemma follows from the fact that the 
map 
$$
	\CIc(G) \to \ind_P^G(\CIc(P)_{\delta})
$$ 
is $R^{\infty}(G)$--linear and the
isomorphism of Shapiro's Lemma,
$$
\cohom_q(G,\ind_P^G(\CIc(P)_{\delta}) \simeq  \cohom_q(P, \CIc(P)_{\delta}),
$$ is natural. 
 
Alternatively, one can use the explicit formula of equation 
\eqref{eq.phiPG}. 
\end{proof}

If $\mfk m =\mfk m_{\gamma}\subset R^{\infty}(G)$ is the maximal  
ideal of functions 
vanishing at a semisimple element $\gamma \in G$, then its image  
$(\rho_P^G)_*(\mfk m):=\rho_P^G(\mfk m)R^{\infty}(P)  
\subset R^{\infty}(P)=R^{\infty}(M)$  
is the ideal of functions 
vanishing at all $g \in M$ that are conjugated to $\gamma$ {\em in} $G$.  
If $\gamma$ is elliptic, then $\mfk m = R^{\infty}(P)$. If $\gamma \in M$, 
then $\mfk m$ need not, in general, be maximal. Nevertheless, we obtain a  
morphism  
$$ 
(\rho_P^G)_{\gamma}:\CC \simeq R^{\infty}(G)_{\gamma}= 
R^{\infty}(G)/{\mfk m}\to  
R^{\infty}(M)/(\rho_P^G)_*(\mfk m) \simeq \CC^{\#\langle\gamma\rangle}, 
$$ 
where $\#\langle\gamma\rangle=l$ is the set 
of conjugacy classes in $M$ that consist 
of elements that are conjugated to $\gamma$ in the bigger group $G$. 
Let $\gamma_1,\gamma_2,\ldots,\gamma_l \in M$ be representatives 
of the conjugacy classes of element in $M$ that are conjugated to  
$\gamma$ in $G$. 
 
We are ready now to study the morphisms 
\begin{multline*}
	(\ind_P^G)_{\gamma}: \Hd_q(\CIc(G))_{\gamma}= \Hd_q(\CIc(G))/
	{\mfk m}\Hd_q(\CIc(G)) \\ \longrightarrow
	\Hd_q(\CIc(P))/(\rho_P^G)_*(\mfk m) \simeq
	\bigoplus\limits_{j=1}^l \Hd_q(\CIc(P))_{\gamma_{j}}.
\end{multline*}

Let $C_P(\gamma_j)$ be the centralizer of $\gamma_j$ in $P$ and 
$C_G(\gamma_j) \simeq C_G(\gamma)$ be the centralizer of $\gamma_j$ 
in $G$. Then $C_P(\gamma_j)_u$ identifies with a subspace of 
$C_G(\gamma)_u$, which gives rise to a continuous proper
map $C_G(\gamma) \times_{C_P(\gamma)} C_P(\gamma)_u 
\to C_G(\gamma)_u$, and hence to a morphism 
$$ 
	\CIc(C_G(\gamma)_u) \to
	\ind_{C_P(\gamma_j)}^{C_G(\gamma)}(\CIc(C_P(\gamma_j)_u)_{\delta})
$$ 
of $C_G(\gamma)$--modules. Passing to cohomology, we obtain using 
Shapiro's Lemma a morphism 
$$ 
	\iota_{\gamma_j}^{\gamma}:
	\cohom_q(C_G(\gamma),\CIc(C_G(\gamma)_u)) \longrightarrow
	\cohom_q(C_P(\gamma),\CIc(C_P(\gamma_j)_u)_{\delta}).
$$ 
 
Recall that Proposition \ref{Prop.local} gives isomorphisms
$$
	\Hd_q(\CIc(G))_{\gamma} \simeq \cohom_q(C_G(\gamma)),
	\CIc(C_G(\gamma)_u))
$$ 
and 
$$
	\Hd_q(\CIc(P))_{\gamma_j} \simeq \cohom_q(C_P(\gamma_j)),
	\CIc(C_P(\gamma_j)_u)_{\delta}).
$$

\begin{proposition}\
Let $\gamma \in G$ be a semisimple element and $M \subset P$ be
as above. If the 
conjugacy class of $\gamma$ does not intersect $M$, then  
$\Hd_*(\CIc(P))_{\gamma}=0$, and hence $(\ind_P^G)_{\gamma}=0$. Otherwise, 
using the notation above, we have 
$$
	(\ind_P^G)_{\gamma} = \oplus_{j=1}^l \iota_{\gamma_j}^\gamma :
	\Hd_*(\CIc(G))_{\gamma} \to \Hd_*(\CIc(P))_{\gamma}\,.
$$ 
\end{proposition}

\begin{proof} This follows from definitions if we observe that in
the sequence of maps
\begin{multline*}
	 G \times_P P \times_{C_P(\gamma_i)} (\gamma_i
	C_P(\gamma_i)_u) \simeq G \times_{C_G(\gamma_i)} C_G(\gamma_i)
	\times_{C_P(\gamma_i)} (\gamma_i C_P(\gamma_i)_u) \\
	\rightarrow G \times_{C_G(\gamma_i)} (\gamma_i
	C_G(\gamma_i)_u)
\end{multline*}
the second map is induced by $C_G(\gamma_i) 
\times_{C_P(\gamma_i)} C_P(\gamma_i)_u \to C_G(\gamma_i)_u$
and their composition induces on homology the direct 
summand $\iota_{\gamma_j}^\gamma$ of the map $(\ind_P^G)_{\gamma}$.
\end{proof}

Another morphism that is likely to play an important role is the
``inflation,'' which we now define. Let $N \subset P$ be the unipotent
radical of an algebraic $p$--adic group, and let $M = P / N$ be its
reductive quotient. Then integration over $N$ defines an algebra
morphism
$$ 
	\psi_M^P:\CIc(P) \to \CIc(M)\, , \; \psi_M^P(f) (m) = \int_N f(mn) dn. 
$$ 
Integration over $N$ also defines a $G$-morphism $\CIc(P)_{\delta} \to
\CIc(M)$, and, since $N$ is a union of compact subgroups, we finally
obtain morphisms
$$ 
\cohom_k(P,\CIc(P)_{\delta}) \to \cohom_k(P,\CIc(M)) \simeq 
\cohom_k(M,\CIc(M)), 
$$ 
whose composition we denote $\inff_M^P$.

\begin{theorem}\ If $M$ is a Levi component of a $p$-adic group
$P$ over a field of characteristic zero, as above. Then we have 
$$
	(\psi_M^P)_*=\inff_M^P : \Hd_*(\CIc(P)) \to \Hd_*(\CIc(M)).
$$
\end{theorem}

\begin{proof}\ Integration over $N$ defines a morphism   
$$
	f:{\mathcal B}(P) \otimes \CIc(P)_{\delta} \to {\mathcal
	B}(M)\otimes \CIc(M),
$$
which commutes with the action of $P$. Then $f \otimes_P 1$ coincides
with the morphism of complexes induced by $\psi_M^P$.

Consider now the maps $h_G$ defined in the proof of Lemma \ref{Lemma.BlBr}.
Then $\psi_M^P \circ h_P = h_M \circ f$, and hence  
$\psi_M^P \circ \tilde{h}_P = \tilde{h}_M \circ (f \otimes_P 1)$,  
from which the result follows.
\end{proof}

We now want to proceed by analogy and establish the explicit form of
the action of $\inff_M^P$ on the geometric fibers. Fix $\gamma \in M$.
Again, integration over the nilpotent radical of $C_P(\gamma)$, the
centralizer of $\gamma$ in $P$, induces a morphism
$$
	\CIc(C_P(\gamma)_u)_{\delta} =
	\CIc(C_P(\gamma)_u)\otimes \Delta_{C_P(\gamma)} \to
	\CIc(C_M(\gamma)_u),
$$
of $P$-modules. Let 
\begin{multline*}
	j_{\gamma} : \Hd_*(\CIc(P))_{\gamma} = \cohom_*(C_P(\gamma),
	\CIc(C_P(\gamma)_u)_{\delta}) \\
	\longrightarrow  \cohom_*(C_M(\gamma),\CIc(C_M(\gamma)_u))
	= \Hd_*(\CIc(M))_{\gamma}
\end{multline*}
be the induced morphism.

\begin{proposition}\ Let $\gamma$ be a semisimple element of a Levi
component $M$ of the group $P$. Let $d(\gamma)$ be the determinant
of $Ad_{\gamma}^{-1} - 1$ acting on $Lie(N)/\ker(Ad_\gamma - 1)$. Then,
using localization at the maximal ideal defined by $\gamma$ in $R^{\infty}(G)=
R^{\infty}(P)$ and the above notation,
$$
	(\inff_M^P)_{\gamma} = |d(\gamma)|^{-1}
	j_{\gamma} : \Hd_*(\CIc(P))_{\gamma} \to
	\Hd_*(\CIc(M))_{\gamma}.
$$
\end{proposition}

\begin{proof}\ Fix $\gamma \in G$, not necessarily
semisimple and let $N_\gamma$ be the subgroup of elements of $N$ commuting
with $\gamma$. We choose a complement $V_\gamma$ of $Lie (N_\gamma)$
in $Lie (N)$ and we use the exponential map to identify 
$V_{\gamma}$ with a subset of $N$.  
Then the Jacobian of the map
$$
	V_{\gamma} \times N_{\gamma} \ni (n,n') \to
	\gamma^{-1} n \gamma n^{-1}n' \in N = V_{\gamma} N_{\gamma}
$$
is $d(\gamma)$, and from this the result follows.
\end{proof} 

This result is compatible with the results of van Dijk on characters
of induced representations, see \cite{vanDijk}.

\section{Examples \label{Sec.Examples}}

Our results can be used to obtain some very explicit results in
certain particular cases.

\vspace*{2mm}{\bf Example 1.}\ 
Let $Z$ be a commutative $p$-adic group of split rank $r$ (so that
$\cohom_q(Z) \simeq \Lambda^q \CC^r$, for all $q \ge 0$).  Then
\begin{equation}
	\Hd_q(\CIc(Z)) \simeq \CIc(Z) \otimes \Lambda^q \CC^r.
\end{equation}

\vspace*{2mm}{\bf Example 2.}\ Let $P$ be the (parabolic) subgroup
of upper triangular matrices in $SL_2(\FF)$, and $A \subset P$ be the
subgroup of diagonal matrices. Then inflation defines a morphism
$$
	\inff_A^P : \Hd_* (\CIc(P)) \to \Hd_*(\CIc(A)) = \CIc(A) \otimes 
	\Lambda^* \CC
$$
whose range is $\CIc(A) \oplus \CIc(A \setminus \{\pm I\})$, with $I$ the
identity matrix of $SL_2(\FF)$. (We see this by localizing at each $\gamma$.)
To describe the kernel of $\inff_A^P$,
let 
\begin{equation}
	u_b = 
\left [ \begin{array}{cc} 1 & b \\ 0 & 1
\end{array}
\right ] .
\end{equation}
Then, if we choose $b$ to range through $\Sigma_u$, a set of representative
of $\FF^*/\FF^{*2}$, the set of elements $u_b$ forms a set of 
representatives of the set of nontrivial conjugacy classes of
unipotent elements of $P$. Recall that $\FF$ has characteristic 
zero, so $\Sigma_u$ is a discrete set. Let $\orb_{b,+}$ be the orbital integral
associated to $u_b$, and let $\orb_{b,-}$ be the orbital integral
associated to $-u_b$, then 
$$
	F_{\pm } = \oplus_{b} \orb_{b,\pm} : 
	\CIc(G) \to \CC^{\Sigma_u}
$$ 
identify the kernel of $\inff_A^P$ as follows.
The map $F_+ \oplus F_{-} : \ker(\inff_A^P) \to \CC^{\pm \Sigma_u}$ is
injective, and the range of each of $F_{\pm}$ is the set of elements
with zero sum.

All in all, let us consider the map $\Phi = \inff_A^P \oplus F_+ 
\oplus F_-$,
$$
	\Phi : \Hd_*(\CIc(P)) \to (\CIc(A) \oplus \CC^{\pm \Sigma_u})_{(0)}
	\oplus (\CIc(A \setminus \{\pm I\}))_{(1)},
$$
where the lower index $(i)$ represents the degree. Then $\Phi$ is 
surjective in degree $1$, and, in degree $0$, it 
consists of $(f, \lambda_{b,\epsilon})$, $f \in \CIc(A)$, 
$\lambda \in \CC$, for $\epsilon \in \{\pm 1\}$ and $b \in 
\Sigma_u \simeq \FF^*/\FF^{*2}$,
such that $\sum_{b} \lambda_{b,\epsilon } = f(\epsilon I)$, for $\epsilon
= \pm 1$. 

Note that evaluation at $\pm I$ does not define a trace on $\CIc(P)$. Actually,
in the spectral sequence of Proposition \ref{Prop.Spectral}, the obstruction
to extend the evaluation at $I$ to a trace is responsible for ``killing''
the $1$ cohomology supported at $I$ (which explains our claim on the
range of $\inff_A^P$ above).

\vspace*{2mm}{\bf Example 3.}\ Consider now the group $G = SL_2(\FF)$,
where $\FF$ is a $p$-adic field of characteristic zero such that
the characteristic of the residual field is not $2$, for
simplicity. Let $F_q$ be the
residual field of $\FF$ (thus $q = p^n$, for some $n \in \NN$ and some prime $p$, 
denotes the number
of elements of $F_q$. We choose $\epsilon$ in the valuation ring
of $\FF$, whose image in $F_q$ is not a square. Also, let $\tau$ be a
generator of the (unique) maximal ideal of the valuation ring of $\FF$.
We shall use the notation of \cite{Sally-Shalika} and thus let $\theta$
range through the set $\{\epsilon, \tau, \epsilon\tau\}$ and let $T_{\theta}$
$T_{\theta}^{\#}$ be the elliptic tori defined there.
(Recall that $T_{\theta}=\{[a_{ij}], a_{11} = a_{22}, a_{21} = 
\theta a_{12}\}$ and $T_{\theta}^{\#} =
\{[a_{ij}], a_{11} = a_{22}, a_{21} = \theta^{\#} a_{12}\}$, where
$\theta^{\#} = \theta a^2$, for some  $a \in \FF^2$ not in the image
of the norm map $N : \FF[ \theta ]^* \to \FF^*$.) We distinguish 
two cases, first the case where $-1$ is a square and then the case where
it is not a square. If $-1$ is a square, then the Weyl group of each
of the tori $T= T_{\theta} , T_{\theta}^{\#}$ has order $2$. Otherwise
$W(T)=\{1\}$, for each $T= T_{\theta}$ or $T= T_{\theta}^{\#}$, but
$T_{\theta}$ and $T_{\theta}^{\#}$ are conjugate for each fixed $\theta$.

Let $X = \cup_{\theta} T_{\theta}/S_2 \cup  \cup_{\theta} 
T_{\theta}^{\#}/S_2$, if $-1$ is a square, and $X = \cup_{\theta}T_{\theta}$\,
otherwise, with the induced topology. Then $X \setminus \{\pm\}$ 
identifies with the set of elliptic conjugacy classes of $SL_2(\FF)$.

Denote by $A \subset SL_2(\FF)$ the set of diagonal matrices in $SL_2(\FF)$.
Let $W(A)=S_2$ act on $\CIc(A) \otimes \Lambda^* \CC$ by conjugation on
$\CIc(A)$ and act by the nontrivial character on $\CC$. Then we
have the following. Recall that there are $10$ conjugacy classes
of unipotent elements in $SL_2(\FF)$, if $p \not = 2$.

\begin{proposition}\ The composition
$$
	\phi:= {\inff_A^P} \circ \ind_P^G : \Hd_*(\CIc(SL_2(\FF))) \to \Hd_*(A)=
	\CIc(A) \otimes \Lambda^* \CC
$$ 
has range consisting of $W(A)$ invariant elements, and the kernel of
$\phi$ is isomorphic to $\CIc(X \setminus \{\pm I\}) \oplus \CC^{10}$,
via orbital integrals with respect to elliptic and unipotent elements.
\end{proposition}

\begin{proof}\ First of all, it is clear that the composition
$\phi = \inff_A^P \circ \ind_P^G$ is invariant with respect to
the Weyl group $W(A)$, and hence its range consists of $W(A)$-invariant
elements. 

The localization of $\phi$ at a regular, diagonal conjugacy class 
$\gamma$ is onto by Proposition \ref{Prop.local}. 
Next, we know that every orbital 
integral extends to $\CIc(SL_2(\FF))$, and this implies directly that
the spectral sequence of Proposition \ref{Prop.Spectral} collapses at the
$E^2$ term. This proves that the localization of $\phi$
at $\gamma = 1$ is also onto, and hence $\phi$ is onto. 
The rest of the proposition follows also from Proposition
\ref{Prop.Spectral} by localization.
\end{proof}  

This example is also discussed in \cite{Blanc-Brylinski}, but from
a different perspective.

\vspace*{2mm}{\bf Example 4.}\ We end this section with a description
of the ingredients entering in the formula \eqref{eq.result} for the
of the Hochschild homology of $\CIc(G)$, if $G=GL_n(\FF)$. Let $\gamma
\in G$ be a semisimple element. The minimal polynomial $Q_{\gamma}$ of
$\gamma$ decomposes as $Q_{\gamma}= p_1 p_2 \ldots p_r$ into
irreducible polynomials with coefficients in $\FF$. (We assume, for
simplicity, that each polynomial $p_j$ is a monic polynomial.) Also,
let $P_{\gamma} = p_1^{l_1} p_2^{l_2} \ldots p_r^{l_r}$ be the
characteristic polynomial of $\gamma$.  Then the algebra generated by
$\gamma$ is
$$
	\FF[\gamma] \simeq \KK_1 \oplus \ldots  \oplus \KK_r,
$$
where $\KK_i = \FF[t]/(p_i(t))$ are not necessarily distinct fields. 
The commutant
$\{\gamma\}'$ of $\gamma$ in $M_n(\FF)$ is the commutant of this
algebra, and hence 
$$
	\{\gamma\}' \simeq M_{l_1}(\KK_1) \oplus M_{l_2}(\KK_2)
	\oplus \dots \oplus M_{l_r}(\KK_r),
$$ 
\begin{equation*}
	C(\gamma) \simeq \prod_{i=1}^r GL_{l_i}(\KK_i) \, , \; \;
	S := Z(C(\gamma)) \simeq \prod_{i=1}^r \KK_i^* , 
\end{equation*}
and
$$
	S^{\rg} = \{(x_i) \in S,\, x_i \text{ generates } \KK_i \text{ and
	the minimal polynomials of }x_i \text{ are distinct} \}.
$$

By the Skolem--Noether theorem, the Weyl group $W(S)= N(S) /C(S)$
coincides with the group of algebra automorphisms of
$\{\gamma\}'$. This group has as quotient a group isomorphic to 
the subgroup $\Pi \subset N(S)$ which permutes the algebras $M_{l_i}(\KK_i)$
Then $\Pi \cong S_{m_1} \times \ldots \times S_{m_t}$, that is, $\Pi$ is a 
product of symmetric groups. We denote the kernel of this morphism by
$W_0(S)$. It is isomorphic to $\prod_{i=1}^r Aut_{\FF}(\KK_i)$ (again by 
the Skolem--Noether theorem). The group $W(S)$ is then the semidirect product 
of $W_0(S)$ by $\Pi$. We hence obtain exact sequences
$$
 	1 \longrightarrow N_0(S) \longrightarrow N(S) \longrightarrow  \Pi 
	\longrightarrow 1 
$$
and
$$
	1 \longrightarrow C(S) \longrightarrow N_0(S) \longrightarrow  W_0(S)
	\longrightarrow 1 \,.
$$

According to \eqref{eq.result}, the only other ingredients necessary to
compute $\Hd_*(\CIc(G))$ are the groups $\cohom_*(C(S), \CIc(\nil
S))$.

Now, the unipotent variety of $C(S)$ is the product of the unipotent
varieties of $GL_{l_i}(\KK_i)$, $i = 1,r$, and the subgroup $C(S)$ preserves
this product decomposition. We see then that in order to prove that the
spectral sequence of Proposition \ref{Prop.Spectral} collapses (for any
choice of open subsets $U_i$), it is
enough to check this for the spectral sequence converging to the cohomology 
of $\CIc(GL_n(\KK)_u)$, for an arbitrary characteristic zero 
$p$-adic field $\KK$.

Fix a unipotent element $\gamma \in GL_n(\KK)$. Define then $V_0 = 0$,
$V_l = \ker (\gamma -1)^l \subset \KK^n$, if $l > 0$. Also, choose
$W_l$ such that $V_l = V_{l-1} \oplus W_l$, and define
$$
	P = \{\gamma \in GL_n(\KK),  \gamma V_l \subset V_l\}, \text{ and }
$$
$$
	M = \{\gamma \in GL_n(\KK),  \gamma W_l = W_l\}.
$$
Then $P$ is a parabolic subgroup with unipotent radical
$$
	N = \{\gamma \in GL_n(\KK),  (\gamma - 1) V_l \subset V_{l-1}\},
$$
and $M$ is a Levi component of $P$. It is easy to check, from
definition, that the $P$--orbit of $u$ in $N$ is dense. The
centralizer of $u$ is then contained in $P$ and has split rank $\le$
the split rank of $P$. Fix a maximal split torus $A$ in the
centralizer of $u$. We can assume that this split torus is contained
in $M$. From the definition and by direct inspection, the map
$\cohom_*(A) \to \cohom_*(M)$ is injective, and hence the map
\begin{equation} \label{eq.surjective}
	\cohom^*(M) \longrightarrow \cohom^*(A) = \cohom^*(C(u))
\end{equation}
is surjective.

Fix now a cohomology class $c_0 \in \cohom^q(C(\gamma)) \simeq
\cohom^q(A)$ and choose a cohomology class $c \in \cohom^q(M)$ that
maps to $c_0$ under the above restriction map. Also, let $\tau$ be the
trace on $\tau_0(f) = f(e)$ on $\CIc(M)$ (obtained by evaluation at
the identity $e$). Then the formula
\begin{equation}
	\phi_0(f_0, \ldots, f_q) = \tau_0(D_c(f_0, \ldots, f_q))
\end{equation}
defines a Hochschild cyclic cocycle on $\CIc(M)$. Consequently,
$$
	\phi = \phi_0 \circ \inff_A^P \circ \circ \ind_P^G
$$
defines a Hochschild cocycle on $\CIc(G)$. For any filtration $U_i$ of
$G_u$ by open, invariant open sets, such that each $U_l \setminus
U_{l-1}$ consists of a single orbit. Suppose that the orbit $U_l
\setminus U_{l-1}$ is the orbit of $\gamma \in GL_n(\FF)$ considered
above.  Then the cocycle $\phi$ will vanish on $\CIc(U_l)$ and
represent the cohomology class 
$$
	c \in \cohom^q(C(\gamma)) \simeq \cohom^q(G, \CIc(U_l
	\setminus U_{l-1})).
$$  
From this it follows that the spectral sequence of 
Proposition~\ref{Prop.Spectral} degenerates at $E^2$.

It is very likely that the above argument extends to arbitrary reductive 
$G$ by choosing $M$ and $P$ as in \cite{Rao}.

\end{document}